\documentclass[preprint]{amsart}

\usepackage{amsmath,amscd,amssymb,amsthm}

\newtheorem{lemma}{Lemma}[section]

\newtheorem{prop}[lemma]{Proposition}
\newtheorem{thm}[lemma]{Theorem}
\newtheorem{cor}[lemma]{Corollary}

\newtheorem{remark}{Remark}[section]

\def\dam{{\da_m}}
\def\dbm{{\db_m}}
\def\dcm{{\dc_m}}
\def\dhm{{\dh_m}}
\def\dlm{{\dl_m}}
\def\lu{{\lambda, \mu}}
\def\ppp{{\phi, \phi'}}
\def\px{\partial_x}
\def\pxi{\partial_{x_i}}
\def\pxj{\partial_{x_j}}

\def\pxr{\partial_{x_r}}
\def\pxe{\partial_\xi}
\def\pxei{\partial_{\xi_i}}
\def\pxej{\partial_{\xi_j}}

\def\pxer{\partial_{\xi_r}}

\def\eval{\mathop{\rm eval}\nolimits}

\def\HC{\mathop{\rm HC}\nolimits}
\def\Lie{\mathop{\rm Lie}\nolimits}
\def\Supp{\mathop{\rm Supp}\nolimits}
\def\tr{\mathop{\rm tr}\nolimits}

\def\gic{{generalized infinitesimal character}}

\def\cal{\mathcal}

\def\D{{\cal D}} 
\def\E{{\cal E}} 
\def\F{{\cal F}}

\def\O{{\cal O}}

\def\S{{\cal S}} 
\def\T{{\cal T}}

\def\frak{\mathfrak}

\def\da{{\frak a}}

\def\db{{\frak b}}

\def\dc{{\frak c}}

\def\dg{{\frak g}}
\def\dgl{\dg\dl}
\def\dh{{\frak h}}

\def\dk{{\frak k}}
\def\dl{{\frak l}}

\def\dn{{\frak n}}

\def\dgp{{\frak p}}
\def\dpgl{\dgp \dgl}

\def\ds{{\frak s}}
\def\dsl{\ds\dl}

\def\dU{{\frak U}}

\def\dZ{{\frak Z}}

\def\Bbb{\mathbb}

\def\bC{\Bbb C}

\def\bN{\Bbb N}

\def\bR{\Bbb R}

\def\bZ{\Bbb Z}

\def\De{\Delta}
\def\ep{\epsilon}

\def\h{\hat}
\def\la{\langle}

\def\o{\overline}
\def\ot{\otimes}
\def\ra{\rangle}

\def\t{\tilde}

\def\Comp{\mathop{\rm Comp}\nolimits}

\def\Div{\mathop{\rm Div}\nolimits}
\def\End{\mathop{\rm End}\nolimits}

\def\Hom{\mathop{\rm Hom}\nolimits}

\def\Mult{\mathop{\rm Mult}\nolimits}

\def\NS{\mathop{\rm NS}\nolimits}
\def\oh{{\ts\frac{1}{2}}}

\def\scrm{\scriptsize\rm}

\def\Span{\mathop{\rm Span}\nolimits}

\def\Symb{\mathop{\rm Symb}\nolimits}
\def\thup{{\mbox{\scrm th}}}
\def\ts{\textstyle}
\def\Vec{\mathop{\rm Vec}\nolimits}
\def\vrm{\Vec \bR^m}

\def\Ext{\mathop{\rm Ext}\nolimits}

\def\bgno{\bigbreak\noindent}
\def\cs{composition series}
\def\dog{differential operator}

\def\hw{highest weight}
\def\hwv{highest weight vector}
\def\ic{{infinitesimal character}}
\def\ie{{\em i.e.,\/}}
\def\iff{if and only if}
\def\ind{indecomposable}
\def\irr{irreducible}
\def\irrep{irreducible representation}

\def\meno{\medbreak\noindent}

\def\pq{projective quantization}
\def\r{representation}

\def\tdm{tensor density module}
\def\tfm{tensor field module}
\def\th{\thinspace}
\def\uea{universal enveloping algebra}
\def\vf{vector field}

\title[Quantization and injective submodules]{Quantization and injective submodules of differential operator modules}

\author{Charles H.\ Conley}
\address{Department of Mathematics \\University of North Texas \\Denton TX 76203, USA} 
\email{conley@unt.edu}

\author{Dimitar Grantcharov}
\address{Department of Mathematics \\University of Texas at Arlington \\Arlington TX 76019, USA}
\email{grandim@uta.edu}

\thanks{The first author was partially supported by Simons Foundation Collaboration Grant 207736.  The second author was partially supported by NSA Grant H98230-13-1-0245.}

\begin{document}
\begin{abstract}
The Lie algebra of \vf s on $\bR^m$ acts naturally on the spaces of \dog s between tensor field modules.  Its projective subalgebra is isomorphic to $\dsl_{m+1}$, and its affine subalgebra is a maximal parabolic subalgebra of the projective subalgebra with Levi factor $\dgl_m$.  We prove two results.  First, we realize explicitly all injective objects of the parabolic category $\O^{\dgl_m}(\dsl_{m+1})$ of $\dgl_m$-finite $\dsl_{m+1}$-modules, as submodules of \dog\ modules.  Second, we study \pq s of \dog\ modules, \ie\ $\dsl_{m+1}$-invariant splittings of their order filtrations.  In the case of modules of \dog s from a tensor density module to an arbitrary tensor field module, we determine when there exists a unique \pq, when there exists no \pq, and when there exist multiple \pq s.

\bgno {\sc 2010 Mathematics Subject Classifications:} 17B10, 17B66
\end{abstract}

\maketitle

\section{Introduction}  \label{Intro} 

It is well known that the \irr\ modules of the classical simple Lie algebras have geometric realizations.  In this paper we give geometric realizations of certain fundamental \ind\ modules: the \ind\ {\em injective\/} objects of the category $\O^{\dgl_m} (\dsl_{m+1})$.  Speaking loosely, this category consists of the $\dgl_m$-finite modules of $\dsl_{m+1}$ whose weights are bounded above.  Its injective modules were first described by Rocha-Caridi \cite{RC80}: either they are parabolic co-Verma modules and have Loewy length~1 or~2, or they are composed of two parabolic co-Verma modules and have Loewy length~3.  The latter type are not co-Verma modules, and to our knowledge they have not previously been realized geometrically.

One natural realization of the \irr\ objects of $\O^{\dgl_m} (\dsl_{m+1})$ is in {\em \tfm s\/} over $\bC^m$.  Indeed, these \irr\ modules are the minimal submodules of the $\dgl_m$-parabolic co-Verma modules, which in turn are those \tfm s coinduced from \irr\ modules of $\dgl_m$.

The \ind\ injective objects of $\O^{\dgl_m} (\dsl_{m+1})$ of Loewy length~3 cannot be realized as \tfm s.  This is because although their highest $\dgl_m$-\r s are semisimple, they do not have \ic s; they have only generalized \ic s.  As shown in Corollary~\ref{injectives are not tfms}, all \tfm s whose highest $\dgl_m$-\r s are semisimple are direct sums of co-Verma modules, and co-Verma modules have \ic s.

In the case that $m = 1$, where $\O^{\dgl_m} (\dsl_{m+1})$ is simply $\O(\dsl_2)$, it has been clear at least since the work \cite{LO99} of Lecomte and Ovsienko that the \ind\ injective objects do have a natural realization: they are contained in modules of linear \dog s between \tfm s.  Our main result, Theorem~\ref{AM Decomp}, proves this statement for all~$m$.

More precisely, we consider the modules of linear \dog s from \tdm s to arbitrary \tfm s, where a {\em \tdm\/} is a \tfm\ coinduced from a scalar module of $\dgl_m$.  In terms of the \hw s of the domain \tdm\ and the range \tfm, we compute explicitly the decomposition of such a module into a direct sum of \ind\ modules of $\dsl_{m+1}$.  We find that every \ind\ injective object of $\O^{\dgl_m} (\dsl_{m+1})$ is realized as a submodule of such a module, and we determine exactly which \dog\ modules contain \ind\ injective objects that are not co-Verma modules.

Our proofs proceed as follows.  In Section~\ref{SC} we compute the action of the Casimir operator of $\dsl_{m+1}$ on the \dog\ modules via the {\em symbol calculus.\/}  The formula we obtain is given in Proposition~\ref{LpppOmegaa} and is closely related to Proposition~7 and Theorem~8 of \cite{MR07}.

In Section~\ref{SMs} we compute the \cs\ of the \dog\ modules by applying the Littlewood-Richardson rule to their symbol modules.  Then we use the theory of the Harish-Chandra homomorphism to prove that generically, no two co-Verma modules appearing in this \cs\ have the same \ic, but in the so-called {\em resonant cases,\/} some \ic s appear twice.  In no case does any \ic\ appear thrice.

In Section~\ref{IOs} we prove that the $\dsl_{m+1}$-intertwining maps between symbol modules may be written explicitly as powers of the divergence operator.  These maps are a particular realization of the $\dsl_{m+1}$-maps between \tfm s, which have been thoroughly analyzed in a more general context by many authors: see \cite{BES88} and the references therein.  The classification of such maps is equivalent to the classification of their duals, the $\dsl_{m+1}$-maps between the $\dgl_m$-finite parabolic Verma modules.  These also have been well studied; see for example \cite{Sh88} and the references therein.

We combine these ingredients in Section~\ref{JORDAN} to show that when there are repeated \ic s, the Casimir operator usually does not act semisimply.  In light of the results of \cite{RC80}, this implies our main result.  As an aside, it also shows that the Casimir operator does not act semisimply on any injective module of Loewy length~3.  We remark that this is not true for all elements of the center of the \uea\ $\dU(\dsl_{m+1})$.  In fact, on any fixed injective module of this type, a subalgebra of the center of codimension~1 acts semisimply.  However, the subalgebra varies with the injective.

There is another application of our analysis of \dog s, in some sense complementary to the realization of injective modules of $\dsl_{m+1}$: it gives conditions for the existence and uniqueness of \pq s.  A {\em \pq\/} of a module of \dog s between two \tfm s is an $\dsl_{m+1}$-invariant splitting of its order filtration, \ie\ an $\dsl_{m+1}$-isomorphism of the module with the direct sum of its symbol modules.

The study of \pq s goes back to \cite{GO96, CMZ97, LO99}.  Most \dog\ modules have unique \pq s.  Those which have either no \pq\ or more than one are said to be {\em resonant.\/}  At first, only \dog s between \tdm s were considered: the {\em scalar case.\/}  This case was completely resolved in \cite{Le00}.

More recently, \pq s of modules of \dog s between arbitrary \tfm s were examined, for example in \cite{Ha07, MR07, CS10}.  In these papers, conditions for resonance are given in terms of the spectrum of the Casimir operator.  The underlying idea is that if a module is resonant, then its Casimir operator has repeated eigenvalues on certain ``tree-like'' submodules; see Proposition~4.10 of \cite{Ha07}, Theorem~11 of \cite{MR07}, and Corollary~6 of \cite{CS10}.

In fact, resonance conditions obtained from the eigenvalues of the Casimir operator are not sharp.  In Remark~\ref{Sec 8 Remark} we show that there are non-resonant modules within which the Casimir operator has repeated eigenvalues on some tree-like submodules.  An elegant abstract necessary and sufficient condition for resonance was obtained in \cite{Mi12}: a module is resonant \iff\ there is a non-trivial $\dsl_{m+1}$-map between two of its symbol modules of different orders.

Our results provide a concrete condition in the case of modules of \dog s from \tdm s to \tfm s: in Proposition~\ref{Resonance Prop} we prove that such a module is resonant \iff\ it has repeated \ic s, and we determine when this occurs.  Theorem~\ref{Multi PQ Thm} resolves the resonant case: the module has no \pq\ \iff\ it contains any \ind\ injective submodules of Loewy length~3.  If it has repeated \ic s but contains no such submodule, then it has multiple \pq s.  Conditions differentiating between the two possibilities are given.

Let us contrast our condition with that of \cite{Mi12}.  For $m>1$, there exist pairs of $\dgl_m$-parabolic co-Verma modules of $\dsl_{m+1}$ with the same \ic, having no $\dsl_{m+1}$-intertwining map between them.  However, the \cs\ of a module of \dog s from a \tdm\ to a \tfm\ never contains such a pair.

Tensor fields and differential operators in fact form modules of the Lie algebra $\vrm$ of polynomial \vf s on $\bR^m$.  This Lie algebra is in one of the four series of infinite dimensional Lie algebras of {\em Cartan type.\/}  Throughout this work we realize $\dsl_{m+1}$ as the {\em projective subalgebra\/} $\dam$ of $\vrm$.  Most of our results concern $\dam$, but parts of Sections~\ref{TFMs & DOs} and~\ref{SC} apply to all of $\vrm$.

The article is organized as follows.  In Section~\ref{SLmpo} we collect the necessary properties of $\dsl_{m+1}$ and the parabolic category $\O^{\dgl_m} (\dsl_{m+1})$, and in Section~\ref{TFMs & DOs} we define the $\vrm$-modules of \tfm s and \dog s.  In Section~\ref{MRs} we state our main results.  All of the remaining sections constitute the proofs of these results.  Section~\ref{SC} reviews the symbol calculus, Section~\ref{SMs} is a study of the symbol modules, and Section~\ref{IOs} reviews $\dsl_{m+1}$-intertwining operators.  Section~\ref{JORDAN} analyzes the Casimir operators of the \dog\ modules, and Section~\ref{PROOFS} concludes the proofs.

\section{Properties of $\dsl_{m+1}$}  \label{SLmpo}

Throughout this article we write $\bN$ for the non-negative integers and $\bZ^+$ for the positive integers.  We establish the notation
\begin{equation*}
   x := (x_1, \ldots, x_m), \qquad
   \bC[x] := \bC[x_1, \ldots, x_m], \qquad
   \pxi := {\partial}/{\partial x_i}.
\end{equation*}
Given a multi-index $I := (I_1, \ldots, I_m) \in \bN^m$, we set
\begin{equation*}
   x^I := x_1^{I_1} \cdots x_m^{I_m}, \qquad
   |I| := I_1 + \cdots + I_m.
\end{equation*}
Given any Lie algebra $\dk$ and any $\dk$-module $V$, the $\dk$-invariant subspace of $V$ is denoted by $V^\dk$:
\begin{equation*}
   V^\dk := \bigl\{ v \in V: X v = 0 \ \forall\ X \in \dk \bigr\}.
\end{equation*}

Let $\vrm$ be the Lie algebra of polynomial \vf s on $\bR^m$, the derivations of the polynomial ring $\bC[x]$:
\begin{equation*}
   \vrm := \Span_\bC \bigl\{ x^I \pxi: I = (I_1, \ldots, I_m) \in \bN^m,\, 1 \le i \le m \bigr\}.
\end{equation*}

We shall always use {\em Einstein's implied summation notation\/} unless explicitly indicated otherwise: in products, repeated indices are summed over from $1$ to $m$.  Thus for example the {\em Euler operator\/} $\E_x := \sum_1^m x_r \pxr$ is written simply $x_r \pxr$.

\subsection{The projective subalgebra}  \label{AM} 

The {\em projective subalgebra\/} $\dam$ of $\vrm$ is
\begin{equation*}
   \dam := \Span_\bC \bigl\{ \pxi,\, x_j \pxi,\, x_j \E_x: 1 \le i, j \le m \bigr\}.
\end{equation*}
It is elementary that $\dam$ is isomorphic to $\dsl_{m+1}$ and $\bigl( \vrm \bigr) \big/ \dam$ is an \irr\ module of $\dam$.  In particular, $\dam$ is a maximal subalgebra of $\vrm$.

We will need several subalgebras of $\dam$:
\begin{equation*} \begin{array}{rcl}
   \dbm &:=& \Span_\bC \bigl\{ \pxi,\, x_j \pxi: 1 \le i, j \le m \bigr\}, \\[6pt]
   \dcm &:=& \Span_\bC \bigl\{ \pxi: 1 \le i \le m \bigr\}, \\[6pt]
   \dlm &:=& \Span_\bC \bigl\{ x_j \pxi: 1 \le i, j \le m \bigr\}, \\[6pt]
   \dhm &:=& \Span_\bC \bigl\{ x_1 \partial_{x_1},\, x_2 \partial_{x_2},
   \ldots,\, x_m \partial_{x_m} \bigr\}.
\end{array} \end{equation*}
Here $\dbm$ is the {\em affine subalgebra,\/} a maximal parabolic subalgebra of $\dam$.  Its Levi factor is $\dlm$, a copy of $\dgl_m$, and its nilradical is the {\em constant subalgebra\/} $\dcm$.  The center of $\dlm$ is spanned by the Euler operator $\E_x$, and $\dh_m$ is a Cartan subalgebra of both $\dlm$ and $\dam$.

In order to give an explicit isomorphism between $\dam$ and $\dsl_{m+1}$, it is convenient to regard $\dsl_{m+1}$ as the quotient $\dpgl_{m+1}$ of $\dgl_{m+1}$ by its center.  Fix coordinates $x_0, \ldots, x_m$ on $\bR^{m+1}$ and write $e_{ij}$ for the elementary matrix with $(i,j)^\thup$ entry~$1$ and all other entries~$0$.  The reader may check that $\phi: \dpgl_{m+1} \to \dam$,
\begin{equation*}
   \phi: e_{00} \mapsto -\E_x,\ e_{0i} \mapsto \pxi,\ e_{i0} \mapsto -x_i \E_x,\ e_{ij} \mapsto x_i \pxj,\ 1 \le i, j \le m.
\end{equation*}
is an isomorphism.  This should be understood geometrically: first check that the map $e_{ij} \mapsto x_i \pxj$ for $0 \le i, j \le m$ is an injective homomorphism from $\dgl_{m+1}$ to $\Vec \bR^{m+1}$.  Then restrict the action of its image to the hyperplane $x_0 = 1$, and identify the hyperplane with $\bR^m$.  Note that $\phi$ carries $\dgl_m$ isomorphically to $\dlm$.

We now describe the root system of $\dam$.  The set $\{ e_{00}, \ldots, e_{mm} \}$ is a basis of a Cartan subalgebra of $\dgl_{m+1}$.  Let $\{ \ep_0, \ldots, \ep_m \}$ be the dual basis.  Then $\{ \ep_1 - \ep_0, \ldots, \ep_m - \ep_0 \}$ is a basis of the dual of the corresponding Cartan subalgebra of $\dpgl_{m+1}$.  Using the isomorphism $\phi$, we regard it as a basis of the dual $\dh_m^*$ of $\dhm$.  Thus, writing $\delta_{ij}$ for the Kronecker delta function,
\begin{equation*}
   \dh_m^* = \Span_\bC \bigl\{ \ep_i - \ep_0: 1 \le i \le m \bigr\}, \quad
   (\ep_i - \ep_0) (x_j \pxj) = \delta_{ij},\ 1 \le i, j \le m.
\end{equation*}

The advantage of this notation is that the action of the Weyl group $W(\da_m)$ of $\dam$ on $\dh_m^*$ is transparent: $W(\da_m)$ is the symmetric group $S_{m+1}$, acting by permutation of the indices $0, 1, \ldots, m$ of the $\ep_i$.  Note that the Weyl group $W(\dl_m)$ of $\dl_m$ is the subgroup $S_m$ of $W(\da_m)$ which permutes the indices $1, \ldots, m$.

The roots of $\da_m$ and $\dl_m$ are
\begin{equation*}
   \Delta(\da_m) = \bigl\{ \ep_i - \ep_j: 0 \le i, j \le m \bigr\} \backslash \bigl\{ 0 \bigr\}, \quad
   \Delta(\dl_m) = \bigl\{ \ep_i - \ep_j: 1 \le i, j \le m \bigr\} \backslash \bigl\{ 0 \bigr\}.
\end{equation*}
The corresponding root vectors are as follows: for $1 \le i, j \le m$,
\begin{equation*}
   \pxi \mbox{\rm\ has root\ } \ep_0 - \ep_i, \qquad
   x_i \E_x \mbox{\rm\ has root\ } \ep_i - \ep_0, \qquad
   x_i \pxj \mbox{\rm\ has root\ } \ep_i - \ep_j.
\end{equation*}

The order $\ep_0 > \ep_1 > \cdots > \ep_m$ fixes a triangular decomposition of $\dam$ and gives the following positive and simple systems.  Below them we list the corresponding positive and simple positive root vectors of $\dam$, and below those we list its negative and simple negative root vectors:
\begin{equation*} \begin{array}{ll}
   \Delta^+(\da_m) = \bigl\{ \ep_i - \ep_j: 0 \le i < j \le m \bigr\}, \quad
   &\Pi^+(\da_m) = \bigl\{ \ep_{i-1} - \ep_i: 1 \le i \le m \bigr\}, \\[6pt]
   \bigl\{ \pxi,\, x_i \pxj: 1 \le i < j \le m \bigr\}, \quad
   &\bigl\{ \partial_{x_1},\, x_1 \partial_{x_2}, \ldots, x_{m-1} \partial_{x_m} \bigr\}, \\[6pt]
   \bigl\{ x_i \E_x,\, x_j \pxi: 1 \le i < j \le m \bigr\}, \quad
   &\bigl\{ x_1 \E_x,\, x_2 \partial_{x_1}, \ldots, x_m \partial_{x_{m-1}} \bigr\}.
\end{array} \end{equation*}

Since $\dl_m$ is a standard Levi subalgebra of $\da_m$, it inherits a triangular decomposition with positive and simple systems
\begin{align*}
   & \Delta^+(\dl_m) = \Delta^+(\da_m) \cap \Delta(\dl_m) = \bigl\{ \ep_i - \ep_j: 1 \le i < j \le m \bigr\}, \\[6pt]
   & \Pi^+(\dl_m) = \Pi^+(\da_m) \cap \Delta(\dl_m) = \bigl\{ \ep_{i-1} - \ep_i: 2 \le i \le m \bigr\}.
\end{align*}

Elements of $\dh_m^*$ are referred to as {\it weights.\/}  We have seen that
\begin{equation*}
   \dh_m^* = \bigl\{ {\ts \sum_0^m} \lambda_i \ep_i:\ {\ts \sum_0^m} \lambda_i = 0 \bigr\}.
\end{equation*}
A weight $\sum_0^m \lambda_i \ep_i$ is {\it $\dam$-dominant\/} if $\lambda_0 \ge \lambda_1 \ge \cdots \ge \lambda_m$, and {\it $\dam$-integral\/} if $\lambda_{i-1} - \lambda_i \in \bZ$ for $1 \le i \le m$.  It is {\it $\dl_m$-dominant\/} if $\lambda_1 \ge \cdots \ge \lambda_m$, and {\it $\dl_m$-integral\/} if $\lambda_{i-1} - \lambda_i \in \bZ$ for $2 \le i \le m$.  (Here $z \ge z'$ for $z, z' \in \bC$ means $z - z' \in \o\bR^+$.)

Given any weight $\lambda$ and any $\dhm$-module $V$, it is conventional to write $V_\lambda$ for the $\lambda$-weight space of $V$.  Given any scalar $c$, it will be convenient to write $V_{(c)}$ for the $c$-eigenspace of the Euler operator $\E_x$.  Thus
\begin{equation*}
   V_\lambda = \bigl\{ v \in V: H v = \lambda(H) v\ \forall\ H \in \dhm \bigr\}, \quad
   V_{(c)} = \bigl\{ v \in V: \E_x v = c v\bigr\}.
\end{equation*}

For $\lambda \in \dh_m^*$, the \irrep\ of $\dlm$ of \hw\ $\lambda$ is denoted by $L_\dlm(\lambda)$.  This module is finite dimensional \iff\ $\lambda$ is $\dlm$-dominant and $\dlm$-integral.  All of the \irrep s of $\dlm$ we will encounter in this article will be finite dimensional.

The algebra $\dcm$ of constant \vf s is $\dlm$-invariant and $\dlm$-\irr: it has \hwv\ $\partial_{x_m}$ and is isomorphic to $L_\dlm(\ep_0 - \ep_m)$.  The \uea\ of $\dcm$ is simply the symmetric algebra $\bC[\px]$ of polynomials in the $\pxi$.  Its homogeneous component of degree~$k$ is its $(-k)$-Euler operator eigenspace, which is $\dlm$-\irr\ and has \hwv\ $\partial_{x_m}^k$.  Thus as $\dlm$-modules,
\begin{equation} \label{Ucm}
   \bC[\px]_{(-k)} = \Span_\bC \bigl\{ \partial_{x_1}^{I_1} \cdots \partial_{x_m}^{I_m}:\ |I| = k \bigr\}
   \simeq L_\dlm(k \ep_0 - k \ep_m).
\end{equation}

\subsection{Infinitesimal characters} \label{ICs}

We now give some information pertaining to the center $\dZ(\dam)$ of $\dU(\dam)$.  
The half-sum $\rho_{\dam}$ of the positive roots of $\dam$ is
\begin{equation*}
   \rho_\dam := \sum_{0 \le i < j \le m} \oh (\ep_i - \ep_j)
   = \oh \bigl( m \ep_0 + (m-2) \ep_1 + \cdots - m \ep_m \bigr).
\end{equation*}
The {\it dot action\/} of the Weyl group $W(\dam)$ on $\dh_m^*$  is
\begin{equation*}
   w \cdot \lambda := w (\lambda + \rho_\dam) - \rho_\dam.
\end{equation*}
Weights with non-trivial stabilizer under the dot action are said to be {\it $\dam$-singular,\/} and all other weights are said to be {\it $\dam$-regular.\/}  A weight $\lambda$ is regular \iff\ $\lambda_0,\, \lambda_1 - 1,\, \ldots,\, \lambda_m - m$ are all distinct.

Regarding elements of $\dU(\dh_m)$ as polynomials on $\dh_m^*$ gives a corresponding dot action on $\dU(\dh_m)$.  We write $\dU(\dhm)^{\cdot W(\dam)}$ for the subspace of dot action invariants.

Let $\dn_m^+$ and $\dn_m^-$ be the direct sums of the positive and negative root spaces of $\dam$, respectively, and note that $\db_m$ contains the Borel subalgebra $\dn_m^+ \oplus \dhm$ of $\dam$.  The zero weight space $\dU(\dam)_0$ of $\dU(\dam)$ is a subalgebra in which $\bigl( \dn_m^- \dU(\dam) \dn_m^+ \bigr)_0$ is a two-sided ideal, and we have the decomposition
\begin{equation*}
   \dU(\dam)_0 = \bigl( \dn_m^- \dU(\dam) \dn_m^+ \bigr)_0 \oplus \dU(\dh_m).
\end{equation*}

The {\it Harish-Chandra homomorphism\/} $\HC: \dU(\dam)_0 \to \dU(\dh_m)$ is the associated projection.   It is a result of Harish-Chandra that $\HC$ restricts to an isomorphism from $\dZ(\dam)$ to $\dU(\dh_m)^{\cdot W(\da_m)}$.  Thus any weight $\lambda \in \dh_m^*$ defines a character
\begin{equation*}
   \chi_\lambda: \dZ(\dam) \to \bC, \qquad
   \chi_\lambda(\Omega) := \HC(\Omega)(\lambda).
\end{equation*}
Characters of $\dZ(\dam)$ are called {\em \ic s\/} of $\dam$.  They are all of the form $\chi_\lambda$, and two characters $\chi_\lambda$ and $\chi_{\lambda'}$ are equal \iff\ $\lambda' \in W(\dam) \cdot \lambda$.

Let $M$ be any $\dam$-module.  We define its {\it $N$-step $\chi_\lambda$-submodule\/} $M^{(\lambda: N)}$, its {\it $\chi_\lambda$-submodule\/} $M^\lambda$, and its {\it generalized $\chi_\lambda$-submodule\/} $M^{(\lambda)}$ by
\begin{equation} \begin{array}{l} \label{gend ICs}
   M^{(\lambda: N)} := \bigl\{ u \in M:
   \bigl( \Omega - \chi_\lambda(\Omega) \bigr)^N u = 0
   \ \forall\ \Omega \in \dZ(\dam) \bigr\}, \\[6pt]
   M^\lambda := M^{(\lambda: 1)}, \qquad
   M^{(\lambda)} := \bigcup_{N=1}^\infty M^{(\lambda: N)}.
\end{array} \end{equation}

If all weight spaces $M_\nu$ of $M$ are finite dimensional and $M$ itself is their direct sum $\bigoplus_{\nu \in \dh_m^*} M_\nu$, then $M$ is a {\it weight module with finite weight multiplicities}.   In this case it is the direct sum of its generalized $\chi_\lambda$-submodules:
\begin{equation} \label{GIC Decomp}
   M = \bigoplus_{\lambda\, \in\, \{ \cdot W(\da_m) \backslash \dh_m^*  \}} M^{(\lambda)}.
\end{equation}

For any semisimple Lie algebra $\dg$, the simplest non-trivial element of $\dZ(\dg)$ is the {\it Casimir element\/} $\Omega_\dg$.  We will need the formulas for $\Omega_\dlm$ and $\Omega_\dam$.  The ordered bases $\{ e_{ij}: 1 \le i, j \le m \}$ and $\{ e_{ji}: 1 \le i, j \le m \}$ of $\dgl_m$ are mutually dual via the trace form, so $\Omega_{\dgl_m} = \sum_{i, j = 1}^m e_{ij} \ot e_{ji}$.  Applying $\phi$ gives
\begin{equation} \label{Omegal}
   \Omega_\dlm = \sum_{i=1}^m \sum_{j=1}^m x_i \pxj \ot x_j \pxi.
\end{equation}
Similarly, $\Omega_{\dgl_{m+1}} = \sum_{i, j = 0}^m e_{ij} \ot e_{ji}$.  Since $\sum_0^m e_{ii} = 0$ in $\dpgl_{m+1}$, we find
\begin{equation*}
   \Omega_{\dpgl_{m+1}} = \Omega_{\dgl_m} +
   \Bigl( \sum_{i=1}^m e_{ii} \Bigr) \ot \Bigl( \sum_{i=1}^m e_{ii} \Bigr)
   + \sum_{i=1}^m \bigl( 2 e_{i,0} \ot e_{0, i} - (m + 1) e_{ii} \bigr).
\end{equation*}
Applying $\phi$ to this gives
\begin{equation} \label{Omegaa}
   \Omega_\dam = \Omega_\dlm + \E_x \ot (\E_x - m - 1)
   - 2\sum_{i=1}^m x_i \E_x \ot \pxi.
\end{equation}
For $\lambda = \sum_{i=0}^m \lambda_i \ep_i$ in $\dh_m^*$, let $||\lambda||^2 := \sum_{i=0}^m \lambda_i^2$.  Standard computations give
\begin{align}
   \HC(\Omega_\dam) &= \E_x \ot (\E_x - m -1) +
   \sum_{i=1}^m x_i \pxi \ot (x_i \pxi + m + 1 - 2i),
   \nonumber \\[6pt] \label{chilambdaOmega}
   \chi_\lambda (\Omega_\dam) &= ||\lambda + \rho_\dam||^2 - ||\rho_\dam||^2.
\end{align}

\subsection{The category $\O$} \label{Obmam}

Here we collect some results on the parabolic category $\O^{\dbm} (\dam)$ of finitely generated $\dbm$-finite weight modules of $\dam$.  For a more detailed treatment of this category, see Chapter~9 of \cite{Hum08} and the references therein.

Henceforth the following notation is in force: $L(\lambda)$ denotes the \irr\ module of $\dam$ of \hw\ $\lambda$, the symbol $\simeq$ denotes equivalence as $\dam$-modules, and if a weight is referred to as dominant, integral, regular, or singular without prefix, it is understood that the prefix is $\dam$.

Given a weight module  $M = \bigoplus_{\nu \in \dh_m^*} M_\nu$ of $\dam$ with finite weight multiplicities, let $M^\vee$ denote the module with underlying space $\bigoplus_{\nu \in \dh_m^*} M_\nu^*$ and action $(X \omega) (v) := \omega(\theta(X) v)$, where $X \in \dam$, $\omega \in M^\vee$, $v \in M$, and  $\theta$ is the Chevalley anti-automorphism  of $\dam$.  Observe that $M$ and $M^\vee$ have the same weights: $(M^\vee)_\nu = M^*_\nu$ because $\theta|_\dhm = 1$.

For $\lambda \in \dh_m^*$, the {\it $\dbm$-parabolic Verma module\/} $M_\dbm(\lambda)$ of $\dam$ is $\dU(\dam) \otimes_{\dU (\dbm)} L_{\dlm}(\lambda)$, where $L_{\dlm}(\lambda)$ is regarded as a $\dbm$-module via the trivial action of the nilradical $\dcm$.  We refer to $M_\dbm^\vee(\lambda)$ as a {\it $\dbm$-parabolic co-Verma module.\/}  It is elementary that $M_{\dbm} (\lambda)$ has a unique irreducible quotient and $M_\dbm^\vee(\lambda)$ has a unique \irr\ submodule, both isomorphic to $L(\lambda)$.  Moreover, both $M_\dbm (\lambda)$ and $M_\dbm^\vee (\lambda)$ have \ic\ $\chi_\lambda$: in the notation of~(\ref{gend ICs}),
\begin{equation} \label{Verma ICs}
   M_\dbm (\lambda) = M_\dbm(\lambda)^\lambda, \qquad
   M_\dbm^\vee (\lambda) = M_\dbm^\vee(\lambda)^\lambda.   
\end{equation}

The modules $M_\dbm(\lambda)$ and $M_\dbm^\vee(\lambda)$ are objects of $\O^{\dbm} (\dam)$ \iff\ $\lambda$ is $\dlm$-dominant and $\dlm$-integral.  In this case, either they are \irr\ or they have Jordan-H\"older length~2.  The next proposition recalls their structure.  The lemma preceding it is a routine exercise in Type~A root systems.

For $\mu \in \dh^*$ and $1 \le i \le m$, define $\mu[0] := \mu$ and $\mu[i] := (0 \cdots i) \cdot \mu$.  Observe that $(0\th i) \cdot \mu[i] = \mu[i-1]$, and if $\mu$ is regular, the $\mu[i]$ are distinct.  Moreover,
\begin{equation} \label{mimommi} \begin{array}{l}
   \mu[i] = (\mu_i - i,\, \mu_0 + 1, \ldots, \mu_{i-1} + 1,\, \mu_{i+1}, \ldots, \mu_m), \\[6pt]
   \mu[i - 1] - \mu[i] = (\mu_{i-1} - \mu_i +1) (\ep_0 - \ep_i).
\end{array} \end{equation}

\begin{lemma} \label{mu bracket i}
Let $\lambda$ be regular, integral, and $\dlm$-dominant.  Then $\lambda = \mu[i]$ for a unique dominant integral weight $\mu$ and a unique $i \in \{0,\ldots, m\}$.
\end{lemma}

\begin{remark} \label{Sec 2.3 Remark}
There are at most $m+1$ $\dlm$-dominant weights in any $W(\dam)$-dot orbit, so there are at most $m+1$ $\dbm$-parabolic Verma modules with any given \ic.  The same is true of\/ $\dbm$-parabolic co-Verma modules.
\end{remark}

\begin{prop} \label{par-Vermas}
Let $\lambda \in \dh_m^*$ be\/ $\dlm$-dominant and\/ $\dlm$-integral.

\begin{enumerate}

\item[(i)]
If $\lambda$ is not both regular and integral, then $M_\dbm(\lambda) \simeq M_\dbm^\vee(\lambda) \simeq L(\lambda)$.

\smallbreak\item[(ii)]
If $\lambda$ is regular integral, set $\lambda = \mu[i]$ for $\mu$ dominant integral and $0 \le i \le m$.

\begin{enumerate}

\smallbreak\item[(a)]
If\/ $0 \le i < m$, then there are non-split exact sequences of\/ $\dam$-modules
\begin{align*}
   & 0 \to L(\mu[i + 1]) \to M_\dbm(\mu[i]) \to L(\mu[i]) \to 0, \\[6pt]
   & 0 \to L(\mu[i]) \to M_\dbm^\vee(\mu[i]) \to L(\mu[i + 1]) \to 0.
\end{align*}
Here $M_\dbm(\mu[i])$ and $M_\dbm^\vee(\mu[i])$ are the only \ind\ modules composed of $L(\mu[i])$ and $L(\mu[i+1])$.

\smallbreak\item[(b)]
If\/ $i = m$, then $M_\dbm(\mu[i]) \simeq M_\dbm^\vee (\mu[i]) \simeq L(\mu[i])$.

\end{enumerate}
\end{enumerate}
\end{prop}

The category $\O^{\dbm} (\dam)$ has enough injectives (and projectives), and we denote by $I(\lambda)$ the injective envelope of $L(\lambda)$ in $\O^{\dbm} (\dam)$.  The following description of the injectives was first established in \cite{RC80}.

\begin{prop} \label{par-injectives}
Let $\lambda \in \dh_m^*$ be\/ $\dlm$-dominant and\/ $\dlm$-integral.

\begin{enumerate}

\item[(i)]
If\/ $\lambda$ is not both regular and integral, then $I(\lambda) \simeq L(\lambda)$.

\smallbreak\item[(ii)]
If\/ $\lambda$ is regular integral, set $\lambda = \mu[i]$ for $\mu$ dominant integral and\/ $0 \le i \le m$.

\begin{enumerate}

\smallbreak\item[(a)]
If\/ $i = 0$, i.e.,\ $\lambda$ is dominant integral, then $I(\lambda) \simeq M_\dbm^\vee (\lambda)$.

\smallbreak\item[(b)]
If\/ $0 < i \le m$, there is a non-split exact sequence of\/ $\dam$-modules
\begin{equation} \label{I no split}
   0 \to M_\dbm^\vee (\mu[i]) \to I(\mu[i]) \to M_\dbm^\vee (\mu[i-1]) \to 0.
\end{equation}

\end{enumerate}
\end{enumerate}
\end{prop}

Just as the non-split exact sequences in Proposition~\ref{par-Vermas}(ii) determine their middle terms, the non-split exact sequence in Proposition~\ref{par-injectives}(ii) determines its middle term.  More precisely, we have the following result concerning extensions between parabolic co-Verma modules.  It is a special case of the formula for the dimensions of all Ext-groups between parabolic Verma modules in the Hermitian symmetric case, deduced in \cite{Sh88}.

\begin{prop} \label{Shelton}
For $\mu$ a dominant integral weight and $0 \le i, j \le m$,
$$
   \dim \Ext^1_{\O^{\dbm}(\dam)}
   \bigl( M_\dbm^\vee (\mu[i]), M_\dbm^\vee (\mu[j]) \bigr) = \delta_{i-1,j}.
$$ 
In particular, if\/ $0 < i \le m$ and 
$$ 
   0 \to M_\dbm^\vee (\mu[i]) \to M \to M_\dbm^\vee (\mu[i-1]) \to 0
$$
is a non-split exact sequence of\/ $\dam$-modules, then $M \simeq I (\mu[i])$.
\end{prop}

Observe that for $\mu$ regular dominant integral and $0 < i < m$, the injective object $I(\mu[i])$ is an \ind\ module of Jordan-H\"older length~4.  In fact, it is a biserial module of Loewy length~3.  On the other hand, $I(\mu[m])$ is a uniserial \ind\ module: its Jordan-H\"older and Loewy lengths are both~3.  In all cases with $i > 0$, $I(\mu[i])$ is projective as well as injective.

Combining~(\ref{gend ICs}), (\ref{Verma ICs}), and Proposition~\ref{par-injectives}, it is clear that $I(\mu[i]) = I(\mu[i])^{(\mu: 2)}$.  However, in general $I(\mu[i])$ does not have an \ic; only a {\em \gic.\/}  Indeed, the following result is ``folklore''.  We thank Volodymyr Mazorchuk for communicating to us a direct proof of Part~(ii).

\begin{prop} \label{folklore}
Assume that $0 < i \le m$.

\begin{enumerate}

\item[(i)]
$I(\mu[i])^\mu$ is the maximal proper submodule of $I(\mu[i])$ and has Loewy length~2.

\smallbreak\item[(ii)]
The Casimir operator $\Omega_\dam$ does not act semisimply on $I(\mu[i])$.

\end{enumerate}
\end{prop}

One should note that the injective-projectives of $\O^{\dbm}(\dam)$ have numerous connections with other fundamental theories, including Khovanov algebras and Brauer tree algebras; see for example Section~4 of \cite{MS11} and the references therein.

\subsection{Tensor products} \label{TPs}

We conclude Section~\ref{SLmpo} with a lemma applicable to any reductive Lie algebra $\dg$.  Although it is undoubtedly known, we have been unable to find a reference, so we give an elementary proof.  Choose a Cartan subalgebra $\dh$ and a positive root system $\Delta^+(\dg)$, and let $\dn^+ \oplus \dh \oplus \dn^-$ be the associated triangular decomposition of $\dg$.

Let $\lambda \in \dh^*$ be any weight and let $V$ be any $\dh$-module.  As is usual, we write $V_\lambda$ for the $\lambda$-weight space of $V$ and $\Supp(V)$ for the set of weights of $V$.  Let $L(\lambda)$ be the \irrep\ of $\dg$ of \hw\ $\lambda$.

Suppose that $\lambda$ and $\mu$ are dominant integral weights, so that $L(\lambda)$ and $L(\mu)$ are finite dimensional.  Observe that
\begin{equation*}
   \bigl( L(\lambda) \ot L(\mu) \bigr)_\nu
   = \bigoplus_{\sigma \in \dh^*} L(\lambda)_\sigma \ot L(\mu)_{\nu - \sigma}.
\end{equation*}
Let $P_{\nu, \sigma} (\lu)$ be the projection to $L(\lambda)_\sigma \ot L(\mu)_{\nu - \sigma}$ along $\bigoplus_{\sigma' \not= \sigma} L(\lambda)_{\sigma'} \ot L(\mu)_{\nu - \sigma'}$:
\begin{equation*}
   P_{\nu, \sigma} (\lu): \bigl( L(\lambda) \ot L(\mu) \bigr)_\nu \twoheadrightarrow
   L(\lambda)_\sigma \ot L(\mu)_{\nu - \sigma}.
\end{equation*}
Then for all $w \in \bigl( L(\lambda) \ot L(\mu) \bigr)_\nu$, we have $w = \sum_\sigma P_{\nu, \sigma} (\lu) (w)$.

\begin{lemma} \label{tensor HWVs}
The restriction of\/ $P_{\nu, \lambda} (\lu)$ to\/ $\bigl( L(\lambda) \ot L(\mu) \bigr)^{\dn^+}_\nu$ is injective.
\end{lemma}

Before the proof let us make some remarks.  Of course, $\bigl( L(\lambda) \ot L(\mu) \bigr)^{\dn^+}_\nu$ is non-zero \iff\ $\nu$ is dominant integral and $L(\nu)$ occurs as a summand of $L(\lambda) \ot L(\mu)$.  A corollary of the lemma is the well-known fact that this can occur only if $\nu \in \lambda + \Supp (L(\mu))$.  The lemma itself may be summarized as ``any \hwv\ of $L(\lambda) \ot L(\mu)$ contains the \hwv\ of $L(\lambda)$ as a factor of a summand''.

\meno {\it Proof.\/}
Suppose that $w$ is a non-zero \hwv\ of weight $\nu$ in $L(\lambda) \ot L(\mu)$.  Let $\sigma_0$ be maximal with respect to the order imposed by $\De^+(\dg)$ such that $P_{\nu, \sigma_0} (w)$ is non-zero.  (We suppress the argument $(\lu)$ of $P_{\nu, \sigma} (\lu)$.)  We will show that $P_{\nu, \sigma_0} (w)$ lies in $L(\lambda)^{\dn^+}_{\sigma_0} \ot L(\mu)_{\nu - \sigma_0}$, which implies $\sigma_0 = \lambda$, proving the result.

For each positive root $\alpha$ in $\De^+(\dg)$, fix a root vector $E_\alpha$ in the $\alpha$-root space $\dn^+_\alpha$.  It will suffice to prove that $P_{\nu, \sigma_0} (w)$ lies in $L(\lambda)^{E_\alpha}_{\sigma_0} \ot L(\mu)_{\nu - \sigma_0}$ for all $\alpha$.

For each $\sigma \in \Supp (L(\lambda))$, let $\bigl\{ v_{\sigma, i}: 1 \le i \le m_\sigma(\lambda) \bigr\}$ be a basis of $L(\lambda)_\sigma$, where $m_\sigma(\lambda)$ is the dimension of $L(\lambda)_\sigma$.  Then there are unique $w_{\sigma, i} \in L(\mu)_{\nu - \sigma}$ such that
\begin{equation*}
   P_{\nu, \sigma} (w) = \sum_{i=1}^{m_\sigma(\lambda)} v_{\sigma, i} \ot w_{\sigma, i}, \qquad
   w = \sum_{\sigma \in \Supp (L(\lambda))}
   \sum_{i=1}^{m_\sigma(\lambda)} v_{\sigma, i} \ot w_{\sigma, i}.
\end{equation*}

By assumption, $E_\alpha w$ and $w_{\sigma_0 + \alpha, i}$ are zero for all $\alpha \in \De^+(\dg)$.  Apply $E_\alpha$ to the above expression for $w$ and use these facts to obtain
\begin{equation*}
   0 = P_{\nu+\alpha, \sigma_0+\alpha} (E_\alpha w) =
   \sum_{i=1}^{m_\sigma(\lambda)} (E_\alpha v_{\sigma_0, i}) \ot w_{\sigma_0, i}.
\end{equation*}

Now fix some choice of $\alpha \in \De^+(\dg)$.  By elementary $\dsl_2$ theory, it is possible to choose the basis $\{ v_{\sigma, i}\}$ of $L(\lambda)$ so that the set of all non-zero vectors of the form $E_\alpha v_{\sigma, i}$ is linearly independent.  For such a choice, the preceding paragraph shows that $w_{\sigma_0, i} = 0$ for all $i$ such that $E_\alpha v_{\sigma_0, i} \not= 0$.  Thus $P_{\nu, \sigma_0} (w)$ lies in $L(\lambda)_{\sigma_0}^{E_\alpha} \ot L(\mu)_{\nu - \sigma_0}$.  Although the basis which gives this result depends on $\alpha$, the result itself holds for all $\alpha$, proving the lemma.  $\Box$

\section{Tensor fields and differential operators}  \label{TFMs & DOs} 

In this section we describe the modules of tensor fields and linear \dog s.  The \tfm s are {\it coinduced\/} from finite dimensional \r s of $\dgl_m$.  When the coinducing module is \irr, the result is a $\dgl_m$-parabolic co-Verma module of $\dsl_{m+1}$.  The co-induction process in fact produces \r s of $\vrm$.  In this paper we are primarily interested in \r s of $\dsl_{m+1}$, but it is natural to give the definitions in full generality.

The module of \dog s between two \tfm s is a special submodule of the entire Hom module between them.  It is of particular interest because it is filtered by order, and the graded symbol modules associated to its filtration are themselves \tfm s.

\subsection{Tensor fields}  \label{TFMs} 

We begin with a concrete realization of \tfm s in spaces of polynomials which is well suited to computation.  Recall that the Einstein summation convention is in force.

\meno {\bf Definition.} {\em
Let $\phi$ be a \r\ of\/ $\dl_m$ on a finite dimensional space $V$.  The corresponding {\em \tfm} of\/ $\vrm$ is the vector space
\begin{equation*}
   \F(V) := \bC[x] \ot V = \bigl\{ h: \bR^m \to V: h \mbox{\it\ a polynomial} \bigr\}.
\end{equation*}
The action of\/ $\vrm$ on $\F(V)$ is the {\em Lie action} $\Lie_\phi$.  It is defined in terms of the\/ {\em $\phi$-divergence,} $\Div_\phi$.  Given a \vf\/ $X = X_j \pxj$ and $h \in \F(V)$,}
\begin{align*}
   & \Div_\phi: \vrm \to \bC[x] \ot \End(V),
   \qquad \Div_\phi(X) := (\pxi X_j) \phi(x_i \pxj), \\[6pt]
   & \Lie_\phi(X) := X + \Div_\phi(X), \qquad
   \Lie_\phi(X) (h) = X_j \pxj h + (\pxi X_j) \phi(x_i \pxj) h.
\end{align*}
{\em If\/ $V$ is 1-dimensional, then $\F(V)$ is known as a {\em \tdm.}}

\medbreak
Let us write $\Lie_\phi|_{\da_m}$ explicitly.  In the computation, it is useful to note that for $f \in \bC[x]$ and $X \in \vrm$, $\Div_\phi(f X) = f \Div_\phi(X) + (X_j \pxi f) \phi(x_i \pxj)$.
\begin{equation} \begin{array}{ll} \label{Lphiam}
   \Lie_\phi(\pxi) = \pxi, &
   \Lie_\phi(x_i \pxj) = x_i \pxj + \phi(x_i \pxj), \\[6pt]
   \Lie_\phi(\E_x) = \E_x + \phi(\E_x), &
   \Lie_\phi(x_i \E_x) = x_i \bigl( \E_x + \phi(\E_x) \bigr) + x_j \phi(x_i \pxj).
\end{array} \end{equation}

We will need the formula for the action of $\Omega_\dam$ on $\F(V)$.  The proof of the following lemma is straightforward from~(\ref{Omegal}), (\ref{Omegaa}), and~(\ref{Lphiam}).  The formulas are not new; for example, they may be deduced from Theorem~8 of \cite{MR07}.

\begin{lemma} \label{LphiOmega}
\begin{enumerate}

\item[(i)]
$\Lie_\phi(\Omega_\dlm) = \phi(\Omega_\dlm) + \E_x (\E_x + m - 1) + 2 \phi(x_i \pxj) x_j \pxi$.

\medbreak \item[(ii)]
$\Lie_\phi(\Omega_\dam) = \phi(\Omega_\dlm) + \phi(\E_x) \bigl( \phi(\E_x) - m - 1 \bigr)$.

\end{enumerate}
\end{lemma}

The following intrinsic characterization of \tfm s will be useful.  Observe that if $M$ is any $\vrm$-module, then $M^\dcm$ is an $\dlm$-module.

\begin{lemma} \label{Intrinsic}
\begin{enumerate}

\item[(i)] $\F(V)^{\dc_m} = V$.

\smallbreak \item[(ii)]
A\/ $\vrm$-module $M$ is isomorphic to a \tfm\ \iff\ it is isomorphic to\/ $\F(M^{\dc_m})$.

\smallbreak \item[(iii)]
If\/ $V'$ is an\/ $\dlm$-submodule of\/ $V$, then\/ $\F(V') = \bC[x] \ot V'$ is a\/ $\vrm$-submodule of\/ $\F(V)$.

\end{enumerate}
\end{lemma}

\meno {\it Proof.\/}
For~(i), $\Lie_\phi(\pxi) = \pxi$ implies that the kernel of the action of $\dcm$ on $\F(V)$ is the space of constant polynomials, $V$.  Part~(ii) is immediate from~(i), and~(iii) is clear from the formula for $\Lie_\phi$.  $\Box$

\medbreak
We now relate \tfm s to co-Verma modules.  Let $\Vec^0 \bR^m$ and $\Vec^0_2 \bR^m$ be the subalgebras of $\vrm$ of \vf s vanishing to first and second order at zero, respectively:
\begin{align*}
   \Vec^0 \bR^m &:= \Span_\bC \bigl\{ x^I \pxi: I \in \bN^m,\, |I| \ge 1,\, 1 \le i \le m \bigr\},
   \\[6pt]
   \Vec^0_2 \bR^m &:= \Span_\bC \bigl\{ x^I \pxi: I \in \bN^m,\, |I| \ge 2,\, 1 \le i \le m \bigr\}.
\end{align*}

Note that $\Vec^0 \bR^m = \dlm \oplus \Vec^0_2 \bR^m$, and $\Vec^0_2 \bR^m$ is an ideal in $\Vec^0 \bR^m$.  Fix a \r\ $V$ of $\dlm$ and extend it trivially to a \r\ of $\Vec^0 \bR^m$.  Define $\Hom_{\Vec^0 \bR^m} \bigl( \dU(\vrm), V \bigr)$ to be the space
\begin{equation*}
   \bigl\{ \kappa: \dU(\vrm) \to V: \kappa(\Theta Y) + \phi(Y) \kappa(\Theta) = 0
   \ \forall\ Y \in \Vec^0 \bR^m \bigr\}.
\end{equation*}
It is a $\vrm$-module under the left action $(X \kappa) (\Theta) := -\kappa (X\Theta)$.  The next lemma shows that as such, it is isomorphic to $\F(V)$.  The proof is elementary.

Given $h \in \F(V)$, define $\eval_0 (h) := h(0)$.  For any Lie algebra $\dg$, let $\Theta \mapsto \Theta^T$ be the anti-automorphism of $\dU(\dg)$ that is $-1$ on $\dg$.  Define
\begin{equation*}
\h h: \dU(\vrm) \to V, \qquad \h h(\Theta) := \eval_0 \circ \Lie_\phi(\Theta^T) h.
\end{equation*}

\begin{lemma} \label{Verma}
\begin{enumerate}

\item[(i)]
$\h h \in \Hom_{\Vec^0 \bR^m} \bigl( \dU(\vrm), V \bigr)$.

\smallbreak \item[(ii)]
$h \mapsto \h h$ is a\/ $\vrm$-isomorphism $\F(V) \to \Hom_{\Vec^0 \bR^m} \bigl( \dU(\vrm), V \bigr)$.

\smallbreak \item[(iii)]
$\big\la h, \Theta \ot \lambda \bigr\ra := \lambda \circ \eval_0 \circ \Lie_\phi(\Theta^T) h$ is a non-degenerate\/ $\vrm$-invariant pairing of\/ $\F(V)$ and\/ $\dU(\vrm) \ot_{\dU(\Vec^0 \bR^m)} V^*$.

\end{enumerate}
\end{lemma}

We abbreviate $\F(L_\dlm(\lambda))$ by $\F(\lambda)$, and write $\Div_\lambda$ and $\Lie_\lambda$ for the associated divergence and Lie action.

\begin{prop} \label{dual Vermas}
Let $\lambda \in \dh_m^*$ be $\dlm$-dominant and $\dlm$-integral.

\begin{enumerate}

\item[(i)]
As $\dam$-modules, $\F(\lambda)$ and $M_\dbm^\vee (\lambda)$ are isomorphic.

\smallbreak \item[(ii)]
$\F(\lambda)$ has \hw\ $\lambda$, \hw\ space $L_\dlm(\lambda)_\lambda$, and \ic\ $\chi_\lambda$.  In particular, $\Lie_\lambda(\Omega_\dam)$ is given by~(\ref{chilambdaOmega}).

\end{enumerate}
\end{prop}

\meno {\it Proof.\/}
Lemma~\ref{Verma} and the definition of $M^\vee$ imply~(i), and~(ii) is a restatement of properties of $M_\dbm^\vee(\lambda)$.  $\Box$

\begin{cor} \label{injectives are not tfms}
Suppose that $M$ is a \tfm\ such that $M^\dcm$ is a weight module.  Then $M$ is a direct sum of\/ $\dbm$-parabolic co-Verma modules.
\end{cor}

\meno {\it Proof.\/}
Combine Lemma~\ref{Intrinsic} and Proposition~\ref{dual Vermas}(i).  $\Box$

\subsection{Scalar-valued differential operators} \label{LieD defn}

We now define the $\vrm$-module of linear \dog s, as well as its symbol module.  For $I \in \bN^m$, let $\px^I$ denote $\partial_{x_1}^{I_1} \cdots \partial_{x_m}^{I_m}$.  The algebra of linear \dog s on $\bC[x]$ with polynomial coefficients is
\begin{equation*}
   \D := \bC[x, \px] = \Span_\bC \bigl\{ x ^I \px^J: I, J \in \bN^m \bigr\}.
\end{equation*}
It is a non-commutative algebra filtered by order: the subspace of order~$\le k$ is
\begin{equation*}
   \D^k := \Span_\bC \bigl\{ x ^I \px^J: I, J \in \bN^m,\, |J| \le k \bigr\}.
\end{equation*}
We write $\Comp$ for composition of \dog s:
\begin{equation*}
   \Comp: \D^{k'} \ot \D^k \to \D^{k' + k}, \qquad
   \Comp(D' \ot D) := D' \circ D.
\end{equation*}

The graded algebra $\S$ associated to the order filtration of $\D$ is the algebra of {\it symbols.\/}  It is conventional to denote the symbol of $\pxi$ by $\xi_i$, so that
\begin{equation*}
   \S := \bC[x, \xi] := \Span_\bC \bigl\{ x ^I \xi^J: I, J \in \bN^m \bigr\}.
\end{equation*}
This algebra is commutative, with $k^\thup$ homogeneous component $\S^k$:
\begin{equation*}
   \S = \bigoplus_{k=0}^\infty \S^k, \qquad
   \S^k := \Span_\bC \bigl\{ x ^I \xi^J: I, J \in \bN^m,\, |J| = k \bigr\}.
\end{equation*}
We write $\Mult$ for symbol multiplication and $\Symb^k$ for the canonical projection:
\begin{equation*} \begin{array}{ll}
   \Mult: \S^{k'} \ot \S^k \to \S^{k' + k}, &
   \Mult(\Xi' \ot \Xi) := \Xi' \Xi \\[6pt]
   \Symb^k: \D^k \to \S^k, &
   \Symb^k \bigl(\sum_{|J| \le k} D_J(x) \px^J \bigr) := \sum_{|J| = k} D_J(x) \xi^J.
\end{array} \end{equation*}
By definition, $\Mult \circ (\Symb^k \ot \Symb^{k'}) = \Symb^{k + k'} \circ \Comp$.  In other words, the symbol of the composition of two operators is the product of their symbols.

The natural action of $\vrm$ on $\D$ is the adjoint action:
\begin{equation*}
   \Lie(X)(D) := X \circ D - D \circ X.
\end{equation*}
The commutativity of symbol multiplication implies that this action preserves the order filtration, and so it drops to an action $\Lie^\S$ of $\vrm$ on $\S$.  The maps $\Comp$, $\Mult$, and $\Symb^k$ are all $\vrm$-covariant.

\begin{remark} \label{Sec 3.2 Remark}
We give explicit formulas for $\Lie$ and $\Lie^\S$ in Section~\ref{SSC}.  In particular, by~(\ref{basic sym powers}), $\S^k$ is a co-Verma module.  Therefore $\D^k$ is in $\O^\dbm(\dam)$, as it is composed of the symbol modules $\S^j$ with $j \le k$.

It is known that $\D^k$ and $\bigoplus_0^k \S^j$ are isomorphic under $\dam$, but for $k > 1$, not under $\vrm$ \cite{LO99}.  Therefore by Lemma~\ref{Intrinsic}(ii), $\D^k$ is not a \tfm\ under $\vrm$, but it is a direct sum of co-Verma modules under $\dam$.
\end{remark}

\subsection{Vector-valued differential operators} \label{LieDppp defn}

In order to define the modules of linear \dog s and symbols between \tfm s, fix two finite dimensional \r s $(\phi, V)$ and $(\phi', V)$ of $\dlm$.  Let $\D(V, V')$ be the space of linear \dog s from $\F(V)$ to $\F(V')$ with polynomial coefficients, and let $\{ \D^k(V, V') \}_k$ be its order filtration:
\begin{equation*}
   \D(V, V') := \D \ot \Hom(V, V'), \qquad
   \D^k(V, V') := \D^k \ot \Hom(V, V').
\end{equation*}
Elements of $\D(V, V')$ are maps from $\F(V)$ to $\F(V')$ as follows: for $D \in \D$, $\tau \in \Hom(V, V')$, $f \in \bC[x]$, and $v \in V$,
\begin{equation*}
   (D \ot \tau) (f \ot v) := D(f) \ot \tau(v).
\end{equation*}

Let $(\phi'', V'')$ be a third \r\ of $\dl_m$.  Composition is
\begin{equation*} \begin{array}{l}
   \Comp: \D^{k'}(V', V'') \ot \D^k(V, V') \to \D^{k' + k}(V, V''), \\[6pt]
   \Comp \bigl( (D' \ot \tau') \ot (D \ot \tau) \bigr) := (D' \circ D) \ot (\tau' \circ \tau).
\end{array} \end{equation*}
Note that for $V = V'$, $\D(V,V)$ is an algebra in which $\D \ot 1$ and $1 \ot \End(V)$ commute.  In other words, $\D(V,V) = \D \ot \End(V)$ in the category of algebras.

The graded space $\S(V, V')$ of symbols associated to $\D(V, V')$ and its $k^\thup$ homogeneous component $\S^k(V, V')$ are
\begin{equation*}
   \S(V, V') := \S \ot \Hom(V, V'), \qquad
   \S^k(V, V') := \S^k \ot \Hom(V, V').
\end{equation*}
Again we have multiplication and projection to symbols:
\begin{equation*} \begin{array}{l}
   \Mult: \S^{k'}(V', V'') \ot \S^k(V, V') \to \S^{k' + k}(V, V''), \\[6pt]
   \Mult \bigl( (\Xi' \ot \tau') \ot (\Xi \ot \tau) \bigr) := (\Xi' \Xi) \ot (\tau' \circ \tau); \\[6pt]
   \Symb^k: \D^k(V, V') \to \S^k(V, V'), \quad
   \Symb^k (D \ot \tau) := \Symb^k(D) \ot \tau.
\end{array} \end{equation*}
Just as in the scalar case, when the domains and ranges are compatible, the symbol of the composition of two operators is the product of their symbols.

The {\it Lie action\/} $\Lie_\ppp$ of $\vrm$ on $\D(V, V')$ is the adjoint action
\begin{equation*}
   \Lie_\ppp(X)(T) := \Lie_{\phi'}(X) \circ T - T \circ \Lie_\phi(X).
\end{equation*}
Because the first order parts of $\Lie_\phi(X)$ and $\Lie_{\phi'}(X)$ are both simply $X$, $\Lie_\ppp$ preserves the order filtration and so drops to an action $\Lie_\ppp^\S$ of $\vrm$ on $\S(V, V')$.  The maps $\Comp$, $\Mult$, and $\Symb^k$ remain $\vrm$-covariant.

\begin{remark} \label{Sec 3.3 Remark}
We give explicit formulas for $\Lie_\ppp$ and $\Lie_\ppp^\S$ in Section~\ref{SSC}.  They show that $\S^k(V, V')$ a direct sum of co-Verma modules, whence $\D^k(V, V')$ is in $\O^\dbm(\dam)$, as it is composed of those $\S^j(V, V')$ with $j \le k$.

Our main results show that in certain cases, $\D^k(V, V') \not\simeq \bigoplus_0^k \S^j(V, V')$ as $\dam$-modules, and so $\D^k(V, V')$ is not a direct sum of co-Verma modules.
\end{remark}

\section{Main results}  \label{MRs}

This section contains the statements of our results.  Their proofs occupy the remainder of the paper.  The chief result is Theorem~\ref{AM Decomp}: all the others follow relatively easily from it.

Fix an $\dlm$-dominant $\dlm$-integral weight $\delta$, so that $\F(\delta) \simeq M^\vee_\dbm(\delta)$ is an object of $\O^\dbm(\dam)$.  Recall from Section~\ref{AM} that $\delta = (\delta_0, \ldots, \delta_m)$, where
\begin{equation*}
   - \delta_0 = \delta_1 + \cdots + \delta_m, \qquad
   \delta_i - \delta_{i+1} \in \bN \mbox{\rm\ \ for\ \ } 1 \le i < m.
\end{equation*}
Here $\delta$ is integral \iff\ $\delta_0 - \delta_i$ is integral for all~$i$.

Our results concern modules of \dog s from \tdm s to \tfm s.  In order to discuss \tdm s precisely, recall that the trace representation of $\dl_m$ on $\bC$ is $\tr: x_i \pxj \mapsto \delta_{ij}$.  The standard divergence $X \mapsto \pxi X_i$ on $\vrm$ is therefore $\Div_{\tr} (X)$.

For $\gamma \in \bC$, write $\bC_\gamma$ for the scalar module of $\dlm$ on which $X$ acts by $\gamma \Div_{\tr}(X)$.  This module has weight $\gamma \sum_{i=1}^m (\ep_i - \ep_0)$.  Such weights are said to be {\em scalar.\/}  The corresponding \tdm\ is $\F(\bC_\gamma)$:
\begin{equation*}
   \bC_\gamma := L_\dlm(-m\gamma, \gamma, \ldots, \gamma), \qquad
   \F(\bC_\gamma) := \F(-m\gamma, \gamma, \ldots, \gamma).
\end{equation*}

We will see in Section~\ref{SC} that  the symbol module of $\D\bigl(\bC_\gamma, \bC_\gamma \ot L_\dlm(\delta)\bigr)$ depends only on $\Hom\bigl((\bC_\gamma, \bC_\gamma \ot L_\dlm(\delta) \bigr) = L_\dlm(\delta)$.  Therefore we use the abbreviations
\begin{equation*}
   \D_\gamma(\delta) := \D\bigl( \bC_\gamma, \bC_\gamma \ot L_\dlm(\delta) \bigr), \qquad
   \S(\delta) := \S\bigl( \bC_\gamma, \bC_\gamma \ot L_\dlm(\delta) \bigr).
\end{equation*}

\subsection{The $\dam$-decomposition of\/ $\D_\gamma(\delta)$}  \label{MR I}

Recall from Section~\ref{ICs} the \ic s $\chi_\mu$ of $\dam$.  We now proceed to state Theorem~\ref{AM Decomp}, the description of the generalized $\chi_\mu$-submodules $\D_\gamma(\delta)^{(\mu)}$ of $\D_\gamma(\delta)$.  By~(\ref{GIC Decomp}), it yields the decomposition of $\D_\gamma(\delta)$ into a direct sum of \ind\ $\dam$-submodules.

Our first proposition gives the decomposition of $\S^k(\delta)$ into a direct sum of co-Verma modules.  It follows easily from Lemma~\ref{Intrinsic} and a particularly simple case of the Littlewood-Richardson rule.  The proof is given in Section~\ref{SMs}.

For each multi-index $K \in \bN^m$, define a corresponding weight $\lambda_K$ by
\begin{equation*}
    \ts \lambda_K := \sum_{i=1}^m K_i (\ep_0 - \ep_i)
    = |K| \ep_0 - \sum_{i=1}^m K_i \ep_i.
 \end{equation*}
We associate to $\delta$ the following sets of multi-indices:
\begin{align*}
    & \kappa(\delta) := \bigl\{ K \in \bN^m:\ K_i \le \delta_i - \delta_{i+1} \th
    \mbox{\rm\ for\ } \th 1 \le i < m \bigr\}, \\[6pt]
    & \kappa(\delta, k) := \bigl\{ K \in \kappa(\delta):\ |K| = k \bigr\}.
\end{align*}
Clearly $\kappa(\delta)$ is the set of all $K$ such that $\delta + \lambda_K$ is $\dlm$-dominant.  For example, $\kappa(\delta, k)$ always contains the multi-index $(0, \ldots, 0, k)$.  If $\delta$ is a scalar weight, then $\kappa(\delta, k)$ contains only this multi-index.

\begin{prop} \label{Decomp of Sk}
$\S^k(\delta) \simeq \bigoplus_{K \in \kappa(\delta, k)} \F(\delta + \lambda_K)$ as $\vrm$-modules.
\end{prop}

Our next proposition gives conditions under which $\S(\delta)$ has repeated \ic s.  Again, the proof is not deep: combining Proposition~\ref{Decomp of Sk} with Harish-Chandra's classification of \ic s, one looks for choices of $K$ and $K'$ such that $\delta + \lambda_K$ and $\delta + \lambda_{K'}$ are in the same dot-orbit of the Weyl group.  The details are given in Section~\ref{SMs}.

In order to state inequalities conveniently, we set $\delta_{m+1} := -\infty$.

\begin{prop} \label{RIC Prop}
Fix $K \in \kappa(\delta)$.  Consider the conditions
\begin{equation} \label{RIC conditions} 
   \delta \mbox{\it\ integral,} \qquad
   \delta + \lambda_K \mbox{\it\ regular,} \qquad
   |K| < \delta_1 - \delta_0.
\end{equation}

\begin{enumerate}

\item[(i)]
If~(\ref{RIC conditions}) does not hold, then $\chi_{\delta + \lambda_{K'}} \not= \chi_{\delta + \lambda_K}$ for all $K' \not= K$ in $\kappa(\delta)$.

\smallbreak \item[(ii)]
If~(\ref{RIC conditions}) holds, then there is a unique multi-index $K' \not= K$ in $\kappa(\delta)$ such that $\chi_{\delta + \lambda_{K'}} = \chi_{\delta + \lambda_K}$, the {\em partner} of $K$.  It is given by
\begin{equation*}
   K'_i = \delta_i - \delta_0 - i - |K|, \qquad K'_j = K_j \mbox{\rm\ \ for\ \ } j \not= i,
\end{equation*}
where $i := i(\delta, |K|)$ is the unique integer in $\{1, \ldots m\}$ such that
\begin{equation*}
   \delta_i - \delta_0 - i \ge |K| > \delta_{i+1} - \delta_0 - (i+1).
\end{equation*}

\smallbreak \item[(iii)]
In~(ii), the dot action of\/ $(0\, i)$ exchanges $\delta + \lambda_K$ and $\delta + \lambda_{K'}$, and
\begin{equation*}
   \lambda_K - \lambda_{K'}
   = (|K| + K_i - \delta_i + \delta_0 + i) (\ep_0 - \ep_i) \not= 0.
\end{equation*}
There is a unique dominant integral weight $\mu$ such that:
\begin{align*}
   & \mbox{\it If\/ $|K| + K_i > \delta_i - \delta_0 - i$,
   then $\delta + \lambda_K = \mu[i-1]$ and
   $\delta + \lambda_{K'} = \mu[i]$.} \\[6pt]
   & \mbox{\it If\/ $|K| + K_i < \delta_i - \delta_0 - i$,
   then $\delta + \lambda_K = \mu[i]$ and
   $\delta + \lambda_{K'} = \mu[i-1]$.}
\end{align*}

\end{enumerate}
\end{prop}

Note that $K$ and $K'$ are mutual partners.  By exchanging them if necessary, we may without loss of generality assume $|K| > |K'|$.

We now state the main result of the paper, the description of the generalized $\chi_\mu$-submodules of $\D_\gamma(\delta)$.  As a corollary, we give the simplest occurrence in a \dog\ module of any given injective object of $\O^\dbm(\dam)$.

\begin{thm} \label{AM Decomp}
Fix $K \in \kappa(\delta)$ and set $k := |K|$.  Let $\simeq$ denote $\dam$-equivalence.
\begin{enumerate}

\item[(i)]
If~(\ref{RIC conditions}) does not hold, then
\begin{equation*}
   \D_\gamma(\delta)^{(\delta + \lambda_K)} =
   \D^k_\gamma(\delta)^{\delta + \lambda_K} \simeq
   \F(\delta + \lambda_K), \qquad
   \D^{k-1}_\gamma(\delta)^{(\delta + \lambda_K)} = 0.
\end{equation*}

\smallbreak \item[(ii)]
If~(\ref{RIC conditions}) holds, set $i = i(\delta, k)$, let $K'$ be the partner of $K$, and set $k' := |K'|$.  Exchanging $K$ and $K'$ if necessary, assume $k >  k'$.  Let $\mu$ be the dominant integral weight such that $\delta + \lambda_K = \mu[i-1]$ and $\delta + \lambda_{K'} = \mu[i]$.

\begin{enumerate}

\smallbreak \item[(a)]
If\/ $k + (m + 1) \gamma \not\in \{1, 2, \ldots, k - k' \}$, then
\begin{equation*}
   \D_\gamma(\delta)^{(\mu)} = 
   \D^k_\gamma(\delta)^{(\mu: 2)} \simeq I(\mu[i]).
\end{equation*}

\smallbreak \item[(b)]
If\/ $k + (m + 1) \gamma \in \{1, 2, \ldots, k - k' \}$, then
\begin{equation*}
   \D_\gamma(\delta)^{(\mu)} = 
   \D^k_\gamma(\delta)^{\mu} \simeq \F(\mu[i-1]) \oplus \F(\mu[i]).
\end{equation*}

\smallbreak \item[(c)]
For all $\gamma$, $\D^{k-1}_\gamma(\delta)^{(\mu)} = \D^{k'}_\gamma(\delta)^\mu \simeq \F(\mu[i])$, and $\D^{k'-1}_\gamma(\delta)^{(\mu)} = 0$.

\smallbreak \item[(d)]
$k - k' = \mu_{i-1} - \mu_i +1$.

\end{enumerate}
\end{enumerate}
\end{thm}

\begin{remark} \label{Sec 4.1 Remark}
There are many choices of\/ $\gamma$, $\delta$, $k$, and $k'$ in this theorem such that $\D^k_\gamma(\delta)^{(\mu)} \simeq I(\mu[i])$, but by~(ii)(d), $k - k'$ depends only on $\mu$ and $i$.
\end{remark}

\begin{cor} \label{main cor}
$I(\delta) \simeq \D_\gamma(\delta)^{(\delta)}$ for all $\gamma \not\in \frac{-1}{m+1} \bN$.  More precisely:

\begin{enumerate}

\item[(i)]
If\/ $\delta$ is either dominant integral or not regular integral, then for all $\gamma$,
\begin{equation*}
   \D_\gamma(\delta)^{(\delta)} = \D_\gamma^0(\delta) \simeq \F(\delta) \simeq I(\delta).
\end{equation*}

\smallbreak \item[(ii)]
If\/ $\delta = \mu[i]$ for some dominant integral $\mu$ and some $1 \le i \le m$, then for all $\gamma \not\in \frac{-1}{m+1} \{ 0, 1, \ldots, \mu_{i-1} - \mu_i \}$,
\begin{equation*}
   \D_\gamma(\delta)^{(\delta)} = \D_\gamma^{\mu_{i-1} - \mu_i +1}(\delta)^{(\delta: 2)} \simeq I(\delta).
\end{equation*}

\end{enumerate}
\end{cor}

\subsection{Projective quantizations}  \label{MR PQ}

Section~\ref{MR I} was concerned with the detection of specific injectives within $\D_\gamma(\delta)$.  In this section we determine whether or not $\D_\gamma(\delta)$ contains {\em any\/} injectives or repeated \ic s, as this determines the existence and uniqueness of \pq s.

\meno {\bf Definition.} {\em
A\/ {\em \pq\/} of\/ $\D(V, V')$ is an $\da_m$-invariant splitting of its order filtration.  For\/ $V$ and $V'$ \irr,\/ $\D(V, V')$ is said to be\/ {\em resonant\/} if it has either no \pq\ or more than one \pq.}

\medbreak
Clearly $\D(V, V')$ has a \pq\ \iff\ the exact sequence
\begin{equation} \label{exact seq}
   0 \to \D^{k-1}(V, V') \to \D^k(V, V') \to \S^k(V, V') \to 0
\end{equation}
is $\dam$-split for all~$k$, and the \pq\ is unique \iff\ the splitting is unique for all~$k$.  Let us state formally the result of Michel mentioned in the introduction.

\begin{thm} \label{JPM Thm} \cite{Mi12}
For $V$ and $V'$ \irr, (\ref{exact seq}) has a unique $\dam$-splitting \iff\ there are no non-trivial $\dam$-maps $\S^k(V, V') \to \bigoplus_{k'<k} \S^{k'}(V, V')$.
\end{thm}

In the case of $\D_\gamma(\delta)$, the next two propositions give explicit conditions for unique splitting of~(\ref{exact seq}) and for resonance.  Note in particular that for $V$ scalar, Michel's condition is equivalent to the absence of repeated \ic s, an \`a priori stronger condition.  We will see in Section~\ref{PROOFS} that both propositions follow from an elementary analysis of Proposition~\ref{RIC Prop} and the well-known classification of intertwining maps between co-Verma modules of $\O^\dbm(\dam)$.

\begin{prop} \label{Splitting Prop}
The following conditions are equivalent:

\begin{enumerate}

\item[(i)]
(\ref{exact seq}) fails to have a unique $\dam$-splitting.

\smallbreak \item[(ii)]
$\S^k(\delta)$ and $\bigoplus_{k' < k} \S^{k'}(\delta)$ share an \ic.

\smallbreak \item[(iii)]
For some $K \in \kappa(\delta, k)$, (\ref{RIC conditions}) holds and the partner $K'$ of $K$ has $|K'| < k$.

\smallbreak \item[(iv)]
$\delta$ is integral, and for some $1 \le i \le m$,
\begin{equation} \label{Splitting conditions}
   \delta_i - \delta_0 - i \ge k > \delta_{i+1} - \delta_0 - i , \qquad
   k > \oh (\delta_i - \delta_0 - i).
\end{equation}
\end{enumerate}
\end{prop}

\begin{prop} \label{Resonance Prop}
The following conditions are equivalent:

\begin{enumerate}

\item[(i)]
$\D_\gamma(\delta)$ is resonant.

\smallbreak \item[(ii)]
$\S(\delta)$ has repeated \ic s.

\smallbreak \item[(iii)]
$\delta$ is integral, and $\delta_1 - \delta_0 > i$ for all $1 \le i \le m$ such that $\delta_i = \delta_1$.

\end{enumerate}
\end{prop}

We now determine which $\D_\gamma(\delta)$ have no \pq s and which have multiple \pq s.  Previously this was known only for $\delta$ scalar \cite{Le00}.  The point is that if $\D_\gamma(\delta)$ is resonant, then it has no \pq\ if it contains at least one injective of Loewy length~3, and it has multiple \pq s if it contains none.

\begin{thm} \label{Multi PQ Thm}
Suppose that~(\ref{exact seq}) fails to have a unique $\dam$-splitting.

\begin{enumerate}

\item[(i)]
(\ref{exact seq}) has no $\dam$-splitting \iff\/ $\D^k_\gamma(\delta)$ contains any injective modules of Loewy length~3 whose symbols are of order~$k$.

\smallbreak \item[(ii)]
If\/ $i(\delta, k) < m$, then~(\ref{exact seq}) has no $\dam$-splitting unless $k + (m+1)\gamma = 1$, when it has a continuous family of $\dam$-splittings.

\smallbreak \item[(iii)]
If\/ $i(\delta, k) = m$, then~(\ref{exact seq}) has no $\dam$-splitting unless the following condition holds, when it has a continuous family of $\dam$-splittings:
\begin{equation} \label{Le00 condn}
   k + (m+1)\gamma \in \bigl\{ 1, \ldots, \max \{1, 2k + m + \delta_0 - \delta_1 \} \bigr\}.
\end{equation}

\end{enumerate}
\end{thm}

For $\delta$ scalar, our (\ref{Splitting conditions}) and (\ref{Le00 condn}) match~(10), (11) of \cite{Le00}, respectively, so our Proposition~\ref{Splitting Prop} and Theorem~\ref{Multi PQ Thm} at $i = m$ match his Theorem~8.5.

\begin{remark} \label{tdelta}
In understanding these results, it is helpful to use the notation
\begin{equation*}
   \t\delta := \delta + \rho_\dam - (\delta_0 + \oh m) (\ep_0 + \cdots + \ep_m), \quad   
   \t\delta_i = \delta_i - \delta_0 - i, \quad \t\delta_{m+1} := -\infty.
\end{equation*}

\begin{itemize}

\item
The $\dlm$-dominance of $\delta$ implies $\t\delta_1 > \t\delta_2 > \cdots > \t\delta_m$.

\smallbreak \item
$\delta$ is integral \iff\/ $\t\delta_i \in \bZ$ for $1 \le i \le m$.

\smallbreak \item
In Proposition~\ref{RIC Prop}, $i(\delta, k)$ is defined by $\t\delta_{i(\delta, k)} \ge k > \t\delta_{i(\delta, k) + 1}$, and~(\ref{RIC conditions}) is equivalent to $\t\delta_1 - |K| \in \bN$ and $|K| + K_{i(\delta, |K|)} \not= \t\delta_{i(\delta, |K|)}$.

\smallbreak \item
In~(\ref{Splitting conditions}), $i = i(\delta, k)$ and\/ $2k > \t\delta_i \ge k > \t\delta_{i+1} + 1$.

\smallbreak \item
Regarding Proposition~\ref{Resonance Prop}, define $i(\delta)$ to be maximal such that\/ $\delta_1 = \delta_2 = \cdots = \delta_{i(\delta)}$.  Then $\D_\gamma(\delta)$ is resonant \iff\ $\t\delta_{i(\delta)} \in \bZ^+$.

\smallbreak \item
$i(\delta) = m$ \iff\ $\delta$ is scalar.
\end{itemize}
\end{remark}

\section{The symbol calculus}  \label{SC}

As noted earlier, the remainder of the paper is devoted to the proofs of the results stated in Section~\ref{MRs}.  In this section we recall the symbol calculus, a technique for computing the Lie action of $\vrm$ on \dog s by transferring it to symbols via the normal order.  The aim is to give the action of the Casimir operator $\Omega_\dam$.  This action was previously computed in \cite{MR07}.  We include a derivation here in order to establish notation and because we give it in a different form.

\subsection{Scalar symbols}  \label{SSC}

This section is a continuation of Section~\ref{LieD defn}.  Recall the modules $\D = \bC[x, \px]$ and $\S = \bC[x, \xi]$ of scalar \dog s and symbols, and their $\vrm$-actions $\Lie$ and $\Lie^\S$.  The {\it normal order total symbol\/} is the linear bijection
\begin{equation*}
   \NS: \D \to \S, \qquad
   \NS(x^I \px^J) := x^I \xi^J.
\end{equation*}
It preserves order and symbols in the following sense:
\begin{equation*}
   \NS^{-1} \bigl( \S^k \bigr) \subseteq \D^k, \qquad
   \Symb^k \circ \NS^{-1} \big|_{\S^k} = 1: \S^k \to \S^k.
\end{equation*}

The {\it normal symbol calculus\/} is a means to give certain types of operators on $\D$ expressions amenable to computation.  It consists in conjugating the operators by $\NS$ to obtain operators on symbol spaces, which turn out to be \dog s in the variables $x$ and $\xi$.  In the literature this conjugation is often suppressed, but for precision we shall usually write it explicitly.

Elements of $\D$ may be regarded as \dog s on $\S$ which commute with constant symbols such as $\xi^J$.  For example, $(x^{I'} \px^{J'}) (x^I \xi^J) : = x^{I'} (\px^{J'} x^I) \xi^J$.  Observe that $\exp(\pxei \ot \pxi)$ is an endomorphism of $\S \ot \S$.  The following lemma gives an elegant expression for the composition operator.  It is well known and simple to verify.

\begin{lemma} \label{NS Comp}
\begin{enumerate}

\item[(i)]
$\NS \circ \Comp \circ (\NS^{-1} \otimes \NS^{-1}) = \Mult \circ \exp (\pxei \ot \pxi)$.

\medbreak \item[(ii)]
In particular, $\NS \bigl( \px^I \circ f(x) \bigr) =
\sum_{J \in \bN^m} \prod_{r=1}^m {I_r \choose J_r} (\px^J f) \xi^{I-J}$.

\end{enumerate}
\end{lemma}

Conjugating the action $\Lie$ by $\NS$ yields an action $\Lie^{\NS}$ of $\vrm$ on $\S$:
\begin{equation*}
   \Lie^{\NS}(X) := \NS \circ \Lie(X) \circ \NS^{-1}.
\end{equation*}

\begin{lemma} \label{LieC}
$\Lie^{\NS}(X) = X - \sum_{|I|>0} \frac{1}{I!} \px^I \bigl( \NS(X) \bigr) \pxe^I$.
\end{lemma}

\meno {\it Proof.\/}
Take $X = X_j \pxj$ and $D = \sum_J D_J(x) \px^J$, so that $\NS(X) = X_j \xi_j$ and $\NS(D) = \sum_J D_J \xi^J$.  Lemma~\ref{NS Comp} gives the following formulas, proving the result.
\begin{align*}
   & \NS(X \circ D) = \sum_J \bigl( X_j D_J \xi_j \xi^J + X(D_J) \xi^J \bigr), \\[6pt]
   & \NS(D \circ X) = \sum_{I, J} {\ts\frac{1}{I!}} \px^I (X_j \xi_j) \pxe^I (D_J \xi^J). \quad \Box
\end{align*}

Clearly the symbol action $\Lie^\S$ is the part of $\Lie^{\NS}$ preserving $\xi$-degree:
\begin{equation*}
   \Lie^\S(X) = X - (\pxi X_j) \xi_j \pxei.
\end{equation*}
Observe that $\Lie^\S$ and $\Lie^{\NS}$ are equal on $\db_m$, and both are the identity on $\dc_m$:
\begin{equation*}
   \Lie^\S \big|_{\db_m} = \Lie^{\NS} \big|_{\db_m}: \pxj \mapsto \pxj, \quad
   x_i \pxj \mapsto x_i \pxj - \xi_j \pxei.
\end{equation*}
Thus $\NS$ is a $\dbm$-isomorphism.  Under both $\Lie^\S$ and $\Lie^{\NS}$ we have $\S^{\dc_m} = \bC[\xi]$, and both actions restrict to the same \r\ $\phi_\xi$ of $\dl_m$ on $\bC[\xi]$.  Let us write $\Div_\xi$ for the divergence $\Div_{\phi_\xi}$ associated to $\phi_\xi$.  Then
\begin{equation*}
   \phi_\xi(x_i \pxj) := -\xi_j \pxei,\ \
   \Div_\xi(X) = -(\pxi X_j) \xi_j \pxei,\ \
   \Lie^\S(X) = X + \Div_\xi (X).
\end{equation*}

By~(\ref{Ucm}), the $\dl_m$-submodule $(\S^k)^{\dc_m} = \bC[\xi]_{(-k)}$ of $\S^{\dc_m}$ is \irr\ with \hw\ $k (\ep_0 - \ep_m)$ and \hwv\ $\xi_m^k$.  Therefore
\begin{equation} \label{basic sym powers}
   \S^k = \bC[x] \ot \bC[\xi]_{(-k)} \simeq \F(k\ep_0 - k \ep_m)
   \mbox{\rm\ as $\vrm$-modules.}
\end{equation}

We now compute $\Lie^\S(\Omega_\dam)$ and $\Lie^{\NS}(\Omega_\dam)$.  Set $\E_\xi := \xi_r \pxer$, and note $\phi_\xi(\E_x) = -\E_\xi$.  Since $\S$ is a \tfm, (\ref{Omegal}) and Lemma~\ref{LphiOmega} give
\begin{equation*}
   \phi_\xi(\Omega_\dlm) = \xi_j \pxei \xi_i \pxej = \E_\xi (\E_\xi + m -1), \qquad
   \Lie^\S(\Omega_\dam) = 2 \E_\xi (\E_\xi + m).
\end{equation*}
Using~(\ref{Lphiam}) and Lemma~\ref{LieC}, we find
\begin{equation*}
   \Lie^\S (x_i \E_x) = x_i (\E_x - \E_\xi) - x_j \xi_j \pxei, \quad
   \Lie^{\NS} (x_i \E_x) = \Lie^\S (x_i \E_x) - \E_\xi \pxei.
\end{equation*}

As noted, the part of $\Lie^{\NS} (\Omega_\dam)$ which preserves $\xi$-degree is $\Lie^\S (\Omega_\dam)$.  Since $\Lie^{\NS}|_{\db_m}$ preserves $\xi$-degree, (\ref{Omegaa}) shows that
\begin{equation} \label{neg1part}
   \Lie^{\NS} (\Omega_\dam) - \Lie^\S (\Omega_\dam) =
   -2 \bigl( \Lie^{\NS} (x_i \E_x) - \Lie^\S (x_i \E_x) \bigr) \pxi =
   2 \E_\xi \pxei \pxi.
\end{equation}

The operator $\pxei \pxi: \S \to \S$ occurring here is the {\em symbol divergence,\/} an operator of $\xi$-degree $-1$ usually denoted by $\Div$.  Thus we have
\begin{equation*}
   \Lie^{\NS} (\Omega_\dam) = 2 \E_\xi (\E_\xi + \Div + m).
\end{equation*}
Keep in mind that the terms in this product do not commute: $[\E_\xi, \Div] = -\Div$.

\begin{remark} \label{Sec 5.1 Remark 2}
The standard divergence $\Div_{\tr}(X)$ defined in Section~\ref{MRs} is simply $\Div \circ \NS(X)$, and so on $\S^1$, $\Div = \Div_{\tr}^{\NS}$.  We will comment on ways to regard the $\phi$-divergence $\Div_\phi$ as an operator on symbols in Lemma~\ref{Remark 5.2}.
\end{remark}

\subsection{Vector symbols}  \label{VSC}

This section is a continuation of Section~\ref{LieDppp defn}.  Recall the modules $\D(V, V')$ and $\S(V, V')$, and their $\vrm$-actions $\Lie_\ppp$ and $\Lie_\ppp^\S$.  Here the normal order total symbol is
\begin{equation*}
   \NS: \D(V, V') \to \S(V,V'), \qquad
   \NS(D \ot \tau) := \NS(D) \ot \tau.
\end{equation*}
As in the scalar case, it preserves order and symbols:
\begin{equation*}
   \NS^{-1} \bigl( \S^k(V, V') \bigr) \subseteq \D^k(V, V'), \qquad
   \Symb^k \circ \NS^{-1}|_{\S^k(V, V')} = 1.
\end{equation*}
Transfer $\Lie_\ppp$ via $\NS$ to the isomorphic \r\ $\Lie^{\NS}_{\phi, \phi'}$ on $\S(V, V')$:
\begin{equation*}
   \Lie^{\NS}_{\phi, \phi'}(X) := \NS \circ \Lie_{\phi, \phi'}(X) \circ \NS^{-1}.
\end{equation*}

Since $\S(V, V') = \S \ot \Hom(V, V')$, we may regard elements of
\begin{equation*}
   \bC[x, \px, \xi, \pxe] \ot \End \bigl( \Hom(V, V') \bigr)
\end{equation*}
as operators on it.  In particular, if $\pi$ is any \r\ of $\dlm$ on $\Hom(V, V')$, then $\Div_\pi(X) = (\pxi X_j) \pi(x_i \pxj)$ acts on $\S(V, V')$.  We will need two such \r s: the adjoint action $\hom(\phi, \phi')$, which we abbreviate to $\hom$, and the right action $\rho_\phi$, defined by $\rho_\phi(Y)(\tau) := -\tau \circ \phi(Y)$.

We now generalize Lemma~\ref{LieC} to the vector-valued case.

\begin{lemma} \label{LieVV}
$\Lie_{\phi, \phi'}^{\NS}(X) = X + \Div_{\hom}(X) + \sum_{|I|>0} \frac{1}{I!} \px^I \bigl( \Div_{\rho_\phi} (X) - \NS(X) \bigr) \pxe^I$.
\end{lemma}

\meno {\it Proof.\/}
Take $X = X_j \pxj$ and $T = \sum_J T_J \px^J$, where $T_J \in \bC[x] \ot \Hom(V, V')$.  Then $\NS(X) = X_j \xi_j$ and $\NS(T) = \sum_J T_J \xi^J$.  Note that
\begin{equation*}
   \Lie_{\phi, \phi'}(X)(T) = \bigl( X \circ T - T \circ X \bigr) + \bigl( \Div_{\phi'}(X) \circ T - T \circ \Div_\phi(X) \bigr).
\end{equation*}

Observe that $\Comp$ factors as $\Comp \ot \Comp$ and $\NS$ factors as $\NS \ot 1$.  It follows that Lemma~\ref{NS Comp} applies just as in the proof of Lemma~\ref{LieC} to give
\begin{equation*}
   \NS(X \circ T) = \sum_J \bigl( X_j T_J \xi_j \xi^J + X(T_J) \xi^J \bigr),\ \
   \NS(T \circ X) = \sum_{I, J} {\ts\frac{1}{I!}} \px^I (X_j \xi_j) \pxe^I (T_J \xi^J).
\end{equation*}
Thus $\NS(T) \mapsto \NS(X \circ T - T \circ X)$ is the operator $X - \sum_{|I|>0} \frac{1}{I!} \px^I \bigl( \NS(X) \bigr) \pxe^I$.

It remains to prove that
\begin{equation*}
   \NS(T) \mapsto \NS \bigl( \Div_{\phi'}(X) \circ T - T \circ \Div_\phi(X) \bigr)
\end{equation*}
is the operator $\Div_{\hom}(X) + \sum_{|I|>0} \frac{1}{I!} \px^I \bigl( \Div_{\rho_\phi} (X) \bigr) \pxe^I$.

Since $\Div_{\phi'}(X) = (\pxi X_j) \phi'(x_i \pxj)$ is an order~$0$ operator on $\S(V')$,
\begin{equation*}
   \NS \bigl( \Div_{\phi'}(X) \circ T \bigr) = (\pxi X_j) \phi'(x_i \pxj) \NS(T).
\end{equation*}
By Lemma~\ref{NS Comp},
\begin{align*}
    \NS \bigl( T \circ \Div_\phi(X) \bigr) & =
   \NS \Bigl( \sum_J \bigl( T_J \circ \phi(x_i \pxj) \bigr) \bigl( \px^J \circ (\pxi X_j) \bigr) \Bigr)
   \\[6pt] & =
   - \sum_{I,J} {\ts\frac{1}{I!}} \bigl( \px^I \pxi X_j \bigr)
   \bigl( \rho_\phi(x_i \pxj) (T_J) \bigr) \bigl( \pxe^I \xi^J \bigr).
\end{align*}
The $I=0$ terms combine with the $\Div_{\phi'}(X)$ terms to give $\Div_{\hom}(X) \NS(T)$, and the $|I| > 0$ terms give the $\Div_{\rho_\phi}$ terms in the lemma.  $\Box$

\medbreak
As in the scalar case, $\Lie^{\NS}_{\phi, \phi'}$ preserves the order filtration, and its $\xi$-degree~$0$ component is the symbol action $\Lie^\S_{\phi, \phi'}$:
\begin{equation*}
   \Lie^\S_{\phi, \phi'}(X) = X + \Div_\xi(X) + \Div_{\hom}(X) = X + \Div_{\phi_\xi \ot \hom}(X).
\end{equation*}
Thus $\Lie^\S_\ppp = \Lie_{\phi_\xi \ot \hom}$, and $\S(V, V')$ is a \tfm:
\begin{equation} \label{SymbVVtfm}
   \S(V, V') = \F\bigl( \S(V, V')^{\dc_m} \bigr) = \F\bigl( \bC[\xi] \ot \Hom(V, V') \bigr).
\end{equation}

At this point we are prepared to give explicit formulas for the restrictions of $\Lie^\S_\ppp$ and $\Lie^{\NS}_\ppp$ to $\da_m$.  As before, they are equal on $\db_m$ and the identity on $\dc_m$:
\begin{equation*}
   \Lie^\S_\ppp \big|_{\db_m} = \Lie^{\NS}_\ppp \big|_{\db_m}: \pxj \mapsto \pxj, \quad
   x_i \pxj \mapsto x_i \pxj - \xi_j \pxei + \hom(x_i \pxj).
\end{equation*}
Use~(\ref{Lphiam}) to obtain $\Lie^\S_\ppp (x_i \E_x)$, and then Lemma~\ref{LieVV} to derive $\Lie^{\NS}_\ppp (x_i \E_x)$:
\begin{equation} \label{LieVVa} \begin{array}{l}
   \Lie^\S_\ppp (x_i \E_x) = x_i \bigl( \E_x - \E_\xi + \hom(\E_x) \bigr) - x_j \xi_j \pxei + x_j \hom(x_i \pxj),
   \\[6pt]
   \Lie^{\NS}_\ppp (x_i \E_x) = \Lie^\S_\ppp (x_i \E_x) - \E_\xi \pxei + \rho_\phi(\E_x) \pxei + \rho_\phi(x_i \pxj) \pxej.
\end{array} \end{equation}

Finally we can deduce the actions of $\Omega_{\da_m}$.  As we mentioned, they were obtained in \cite{MR07} in a different form; see Section~3.6 of that paper.  In our applications $\phi$ will be scalar, so we give a simplified expression in that case.  Recall the symbol divergence $\Div = \pxei \pxi$ and the $\dlm$-module $\bC_\gamma$.

\begin{prop} \label{LpppOmegaa}
\begin{enumerate}

\item[(i)]
The symbol action of\/ $\Omega_{\da_m}$ on\/ $\S(V, V')$ is
\begin{align*}
   \Lie^\S_\ppp (\Omega_{\da_m}) = &
   \hom(\Omega_{\dl_m}) - 2 \hom(x_i \pxj) \xi_i\pxej \\[6pt]
   & + \hom(\E_x) \bigl( \hom(\E_x) - 2 \E_\xi - m - 1 \bigr)
   + 2 \E_\xi \bigl( \E_\xi + m \bigr).
\end{align*}

\smallbreak \item[(ii)]
The \dog\ action of\/ $\Omega_{\da_m}$ on\/ $\S(V, V')$ is
\begin{equation*}
   \Lie^{\NS}_\ppp (\Omega_{\da_m}) =
   \Lie^\S_\ppp (\Omega_{\da_m})
   - 2 \rho_\phi(x_i \pxj) \pxej \pxi
   - 2 \bigl( \rho_\phi(\E_x) - \E_\xi \bigr) \Div.
\end{equation*}

\smallbreak \item[(iii)]
If\/ $V = \bC_\gamma$, then\/ $\Lie^{\NS}_\ppp (\Omega_{\da_m}) = \Lie^\S_\ppp (\Omega_{\da_m}) + 2 \bigl( \E_\xi + (m+1) \gamma \bigr) \Div$.
\end{enumerate}
\end{prop}

\meno {\it Proof.\/}
For~(i), apply Lemma~\ref{LphiOmega} to $\Lie_{\phi_\xi \ot \hom}$.  The result follows from the formulas $(\phi_\xi \ot \hom) (\E_x) = \hom(\E_x) - \E_\xi$ and
\begin{align*}
   (\phi_\xi \ot \hom) (\Omega_{\dl_m}) &=
   \bigl( \hom(x_i \pxj) - \xi_j \pxei \bigr) \bigl( \hom(x_j \pxi) - \xi_i \pxej \bigr) \\[6pt]
   &= \hom(\Omega_{\dl_m}) - 2 \hom(x_i \pxj) \xi_i\pxej + \E_\xi (\E_\xi + m - 1).
\end{align*}

For~(ii), the idea leading to~(\ref{neg1part}) gives
\begin{equation*}
   \Lie^{\NS}_\ppp (\Omega_\dam) - \Lie^\S_\ppp (\Omega_\dam) =
   -2 \bigl( \Lie^{\NS}_\ppp (x_i \E_x) - \Lie^\S_\ppp (x_i \E_x) \bigr) \pxi.
\end{equation*}
Now~(\ref{LieVVa}) gives the result.  Part~(iii) is immediate.  $\Box$

\section{The symbol modules}  \label{SMs}

The purpose of this section is to prove Propositions~\ref{Decomp of Sk} and~\ref{RIC Prop}.  Combining~(\ref{basic sym powers}) and~(\ref{SymbVVtfm}), we see that as $\vrm$-modules,
\begin{equation*}
   \S^k(\delta) = \F\bigl( \bC[\xi]_{(-k)} \ot L_\dlm(\delta) \bigr)
   \simeq \F\bigl( L_\dlm(k \ep_0 - k \ep_m) \ot L_\dlm(\delta) \bigr).
\end{equation*}
Recall $\lambda_K$ and $\kappa(\delta, k)$.  The decomposition of $L_\dlm (k \ep_0 - k \ep_m) \ot L_\dlm(\delta)$ into $\dlm$-\irr\ summands is a special case of the Littlewood-Richardson rule:

\begin{prop} \label{LR}
$L_\dlm (\lambda_{k e_m}) \ot L_\dlm(\delta) \th \simeq \th \bigoplus_{K \in \kappa(\delta, k)} L_\dlm(\delta + \lambda_K)$.
\end{prop}

For future reference, at this point we specify a \hwv\ $v_\delta^K$ of weight $\delta + \lambda_K$ for each $K$ in $\kappa(\delta, k)$.  Fix a \hwv\ $v_\delta$ of $L_\dlm(\delta)$.  

\begin{lemma} \label{vdeltaK}
For $K \in \kappa(\delta, k)$, there exists a unique \hwv\ $v_\delta^K$ in\/ $\bC[\xi]_{(-k)} \ot L_\dlm(\delta)$ of weight $\delta + \lambda_K$ with the following properties:
\begin{enumerate}

\item[(i)]
$v_\delta^K = \sum_{|J| = k} \xi^J \ot v_\delta^K(J)$, where $v_\delta^K(J) \in L_\dlm(\delta)_{\delta + \lambda_K - \lambda_J}$.

\smallbreak \item[(ii)]
$v_\delta^K(K) = v_\delta$.

\smallbreak \item[(iii)]
$v_\delta^K(J) = 0$ if\/ $\sum_{i = 1}^{i_0} (K_i - J_i) < 0$ for any~$i_0$.
\end{enumerate}
\end{lemma}

\meno {\it Proof.\/}
Since $\bC[\xi]_{\lambda_J} = \bC \xi^J$, any vector of weight $\delta + \lambda_K$ has a unique expression as in~(i).  We claim that there exists a unique \hwv\ $v_\delta^K$ of weight $\delta + \lambda_K$ satisfying~(ii).  Write $\mu$ for $k (\ep_0 - \ep_m)$, regard $\bC[\xi]_{(-k)} \ot L_\dlm(\delta)$ as $L_\dlm(\delta) \ot L_\dlm(\mu)$, and apply Lemma~\ref{tensor HWVs} to see that the projection
\begin{equation*}
   P_{\delta + \lambda_K, \delta} (\delta, \mu):
   \bigl( L_\dlm(\delta) \ot L_\dlm(\mu) \bigr)_{\delta + \lambda_K} \twoheadrightarrow
   L_\dlm(\delta)_\delta \ot L_\dlm(\mu)_{\lambda_K}
\end{equation*}
is injective on \hwv s.  By Proposition~\ref{LR}, its domain contains a single line of \hwv s, and since $L_\dlm(\mu)$ has 1-dimensional weight spaces, its codomain is the line $\bC v_\delta \ot \xi^K$.  Let $v_\delta^K$ be the unique \hwv\ projecting to $v_\delta \ot \xi^K$.  This proves the claim.

Finally, $L_\dlm(\delta)_{\delta + \lambda_K - \lambda_J}$ is zero unless $\lambda_J - \lambda_K$ is a non-negative integer combination of the simple root vectors $\ep_{i-1} - \ep_i$, so $v_\delta^K$ also satisfies~(iii).  $\Box$

\meno {\it Proof of Proposition~\ref{Decomp of Sk}.\/}
For $K \in \kappa(\delta, k)$, write $L_\dlm(v_\delta^K)$ for the copy of $L_\dlm(\delta + \lambda_K)$ in $\bC[\xi]_{(-k)} \ot L_\dlm(\delta)$ with \hwv\ $v_\delta^K$, and $\F(v_\delta^K)$ for the $L_\dlm(v_\delta^K)$-valued polynomials:
\begin{equation*}
   L_\dlm(v_\delta^K) := \dU(\dlm) v_\delta^K, \qquad
   \F(v_\delta^K) := \bC[x] \ot L_\dlm(v_\delta^K).
\end{equation*}

By Lemma~\ref{Intrinsic}, $\F(v_\delta^K)$ is a $\vrm$-submodule of $\S^k(\delta)$ isomorphic to the \tfm\ $\F(\delta + \lambda_K)$.  To complete the proof, use Proposition~\ref{LR} to obtain
\begin{equation} \label{decomp of Sk}
   \bC[\xi]_{(-k)} \ot L_\dlm(\delta) = \bigoplus_{K \in \kappa(\delta, k)} L_\dlm(v_\delta^K),
   \quad
   \S^k(\delta) = \bigoplus_{K \in \kappa(\delta, k)} \F(v_\delta^K). \ \ \Box
\end{equation}

\meno {\it Proof of Proposition~\ref{RIC Prop}.\/}
Recall $\t\delta$ from Remark~\ref{tdelta}.  As described in Section~\ref{ICs}, $\chi_{\delta + \lambda_K} = \chi_{\delta + \lambda_{K'}}$ \iff\ $\delta + \lambda_K + \rho_\dam$ and $\delta + \lambda_{K'} + \rho_\dam$ are in the same $S_{m+1}$-orbit.  Subtracting the $S_{m+1}$-stable vector $(\delta_0 + \oh m)(1, \ldots, 1)$, the condition is $\t\delta + \lambda_{K'} = w (\t\delta + \lambda_K)$ for some $w$ in $S_{m+1}$.

From the definition of $\kappa(\delta, k)$, we find that
\begin{equation*}
   \t\delta_1 \,\ge\, \t\delta_1 - K_1 \,>\, \t\delta_2 \,\ge\, \t\delta_2 - K_2
   \,>\, \t\delta_3 \,\ge\, \cdots \,>\, \t\delta_m \,\ge\, \t\delta_m - K_m,
\end{equation*}
and similarly for $K'$.  Therefore $\t\delta_i - K'_i \not= \t\delta_j - K_j$ for $i \not= j$, so if $w$ exists, it must be a transposition $(0 \th i)$.  Because $(\t\delta + \lambda_K)_0 = |K|$, this would force
\begin{equation*}
   |K| = \t\delta_i - K'_i, \quad
   |K'| = \t\delta_i - K_i, \quad
   K'_j = K_j \th \mbox{\rm\ for\ } \th j \not= i.
\end{equation*}
These equations imply $\delta$ integral and $\t\delta_i \ge |K| > \t\delta_{i+1}$.  The regularity of $\delta + \lambda_K$ is equivalent to both $K \not= K'$ and $\t\delta_i \not= |K| + K_i$.  Thus if $w$ exists, (\ref{RIC conditions}) holds, and~(ii) is proven.  For~(iii), apply~(\ref{mimommi}) and Lemma~\ref{mu bracket i}.  $\Box$

\section{Intertwining operators}  \label{IOs} 

Here we describe the affine ($\dbm$-covariant) and projective ($\dam$-covariant) operators on \tfm s, and in particular, on symbol modules.

\subsection{Affine operators} \label{bmaps}

Recall the algebra $\bC[\px]$ from~(\ref{Ucm}), and write $\bC[[\px]]$ for the corresponding algebra of formal power series.  Note that $\bC[[\px]]$ acts naturally on $\bC[x]$.

\begin{lemma} \label{cTFMops}
For any\/ $\dl_m$-modules $V$ and $V'$, there is an isomorphism
\begin{equation*}
   \Hom_{\dc_m} \bigl( \F(V), \F(V') \bigr) \simeq \bC[[\px]] \ot \Hom(V, V').
\end{equation*}
It intertwines the $\dl_m$-actions, and for $V = V'$ it is an algebra isomorphism.
\end{lemma}

\meno {\it Proof.\/}
Recall that $\dcm$ acts solely on the first factor of $\F(V) = \bC[x] \ot V$, so
\begin{equation*}
   \Hom_{\dc_m} \bigl( \F(V), \F(V') \bigr) \simeq \End_\dcm \bigl(\bC[x]\bigr) \ot \Hom(V, V').
\end{equation*}
It is elementary that $\End_\dcm (\bC[x]) = \bC[[\px]]$.  For the $\dlm$-actions, use the fact that $\Div_\phi(X) = \phi(X)$ for $X \in \dlm$.  The rest is easy.  $\Box$

\begin{lemma} \label{bTFMops}
Suppose that\/ $\E_x$ acts on\/ $V$ and\/ $V'$ by scalars:\/ $V = V_{(c)}$ and $V' = V'_{(c')}$.  Then\/ $\Hom_\dbm \bigl( \F(V), \F(V') \bigr)$ is~$0$ unless $c' - c \in \bN$, when it is isomorphic to
\begin{equation*}
   \Bigl[ \bC[\px]_{(c-c')} \ot \Hom \bigl( V, V' \bigr) \Bigr]^{\dl_m}.
\end{equation*}
\end{lemma}

\meno {\it Proof.\/}
Use $\Hom(V, V') = \Hom(V, V')_{(c' - c)}$ together with Lemma~\ref{cTFMops}.  $\Box$

\medbreak
By~(\ref{SymbVVtfm}), $\S^k(V, V') = \F \bigl(\bC[\xi]_{(-k)} \ot \Hom(V, V') \bigr)$.  If $\E_x$ acts by scalars on $V$ and $V'$, then it acts by zero on $\End \bigl( \Hom(V, V') \bigr)$.  In this case Lemma~\ref{bTFMops} gives the following corollary describing the $\db_m$-maps between symbol modules.

\begin{cor} \label{bSymbolOps}
Assume that\/ $\E_x$ acts by scalars on $V$ and $V'$.  Then for $k < k'$ the space\/ $\Hom_\dbm \bigl(\S^k(V, V'), \S^{k'}(V, V') \bigr)$ is\/~$0$, and for $k \ge k'$ it is isomorphic to
\begin{equation} \label{Divr}
   \Bigl[ \bC[\px]_{(k' - k)} \ot \Hom\bigl( \bC[\xi]_{(-k)}, \bC[\xi]_{(-k')} \bigr)
   \ot \End \bigl( \Hom(V, V') \bigr) \Bigr]^{\dl_m}.
\end{equation}
\end{cor}

\begin{remark} \label{Sec 7.1 Remark 1}
This corollary describes the associative algebra $\End_\dbm \bigl(\S(V, V') \bigr)$.  For example, if both $V$ and $V'$ are scalar modules, elementary \r\ theory shows that~(\ref{Divr}) is 1-dimensional, so it is spanned by $\Div^{k - k'}$.  This leads to the well-known fact that in this scalar case, the $\dbm$-endomorphism algebra of the symbols is $\bC[\E_\xi, \Div]$: see Section~4.1 of \cite{LO99}.
\end{remark}

Let us elaborate on Remark~\ref{Sec 5.1 Remark 2}.  Suppose that $\pi$ is any \r\ of $\dlm$ on $\Hom(V, V')$ with the property that $[\hom(X), \pi(Y)] = \pi([X,Y])$ for all $X$ and $Y$ in $\dlm$.  Then there is a corresponding $\db_m$-invariant operator $\Div^{\NS}_\pi$ on $\S(V, V')$, the {\it $\pi$-symbol divergence.\/}  It is defined in the following lemma, which is simple to prove directly.

\begin{lemma} \label{Remark 5.2}
$\Div_\pi^{\NS} := \pi(x_i \pxj) \pxej \pxi$ is in\/ $\Hom_{\db_m} \bigl(\S^k(V, V'), \S^{k-1}(V, V') \bigr)$.
\end{lemma}

The notation is chosen to reflect the fact that $\Div_\pi^{\NS} \bigl( \NS(X) \bigr) = \Div_\pi(X)$.  Some examples of such \r s $\pi$ are $\hom$ itself, the left and right actions $\lambda_{\phi'}$ and $\rho_\phi$, and if properly interpreted, any \irr\ constituent of $\hom$.

Similarly, one finds that the operator $\Div_{\phi_\xi}^{\NS} := \phi_\xi(x_i \pxj) \pxej \pxi$ is $\db_m$-invariant.  It simplifies to $-\E_\xi \Div$, so Proposition~\ref{LpppOmegaa}(ii) may be written as
\begin{equation*}
   \Lie^{\NS}_\ppp (\Omega_{\da_m}) - \Lie^\S_\ppp (\Omega_{\da_m}) =
   - 2 \bigl( \Div^{\NS}_{\rho_\phi} + \Div^{\NS}_{\phi_\xi} + \rho_\phi(\E_x) \Div \bigr).
\end{equation*}
This is a $\dbm$-map, as follows from the fact that $\Lie^{\NS}_\ppp|_\dbm$ and $\Lie^\S_\ppp|_\dbm$ are equal and hence commute with both $\Lie^{\NS}_\ppp (\Omega_\dam)$ and $\Lie^\S_\ppp (\Omega_\dam)$.

We now prove a key proposition about the image of $\S(\delta)$ under powers of $\Div$.  We begin with two lemmas.  The first is elementary:

\begin{lemma} \label{Divd}
\begin{enumerate}
\item[(i)]
$\Div^d \circ (x_{i_1} \cdots x_{i_d})$ acts on\/ $\bC[\xi] \ot L_\dlm(\delta)$ as $d! \th \partial_{\xi_{i_1}} \cdots \partial_{\xi_{i_d}}$.

\smallbreak \item[(ii)]
$\Div^d \circ \bigl( \bC[x]_{(d)} \bigr)$ acts on\/ $\bC[\xi] \ot L_\dlm(\delta)$ as\/ $\bC[\partial_\xi]_{(d)}$.

\end{enumerate}
\end{lemma}

The second lemma gives a formula for the \hwv s of the tensor product of the standard \r\ $\bC^m$ of $\dgl_m$ with any finite dimensional \irrep\ of $\dgl_m$.  It is known; we outline two proofs.

Recall from Section~\ref{AM} that $e_{ij}$ is the $ij^\thup$ elementary matrix in $\dgl_m$.  The standard basis of the weight space of $\dgl_m$ is $\ep_1, \ldots, \ep_m$, where the weight $\ep_i$ of $\dgl_m$ is identified with the weight $\ep_i - \ep_0$ of $\dlm$.  It will be convenient to have a distinct notation for the standard basis of $\bC^m$, so we use $e_1, \ldots, e_m$:
\begin{equation*}
   \bC^m = \Span \bigl\{ e_1, e_2, \ldots, e_m \bigr\}.
\end{equation*}   

Given a dominant integral weight $\nu = \sum_1^m \nu_i \ep_i$ of $\dgl_m$ and integers $1 \le j < i \le m$, $r \ge 0$, and $1 \le i_r < \cdots < i_0 \le m$, define scalars
\begin{equation*}
   c_j(i) := \frac{1}{\nu_j - \nu_i + i - j} \,, \quad
   b(i_0, i_1, \ldots, i_r) := (-1)^r \, \prod_{s=1}^r c_{i_s}(i_0) \, \prod_{j=1}^{i_r - 1} \bigl( 1 - c_j(i_0) \bigr).
\end{equation*}

\begin{lemma} \label{HWV}
Fix a \hwv\ $v_\nu$ of\/ $L_{\dgl_m}(\nu)$.  If\/ $\nu_{i_0 - 1} > \nu_{i_0}$, then\/ $\bC^m \ot L_{\dgl_m}(\nu)$ has a \hwv\ $v_\nu^{i_0}$ of weight $\nu + \ep_{i_0}$, given by
\begin{equation*}
   v_\nu^{i_0} := 
   \sum_{0 \le r < i_0} \ \sum_{1 \le i_r < \cdots < i_1 < i_0} \,
   b(i_0, \ldots, i_r) \, e_{i_r} \ot \bigl( e_{i_0 i_1} e_{i_1 i_2} \cdots e_{i_{r-1} i_r} v_\nu \bigr).
\end{equation*}
\end{lemma}

\meno {\it Proof.\/}
This result can be proven using the non-commutative finite factorization of the {\em extremal projector\/} discovered in \cite{AST79}.  Their theorem implies that
\begin{equation*}
   \bigl(1 - c_{i_0 - 1}(i_0) e_{i_0, i_0 - 1} e_{i_0 - 1, i_0} \bigr) \cdots
   \bigl(1 - c_2(i_0) e_{i_0, 2} e_{2, i_0} \bigr) \bigl(1 - c_1(i_0) e_{i_0, 1} e_{1, i_0} \bigr)
   \bigl(e_{i_0} \ot v_\nu \bigr)
\end{equation*}
is a \hwv.  A delicate computation shows that it is $v_\nu^{i_0}$.

For a ``low technology'' proof, first use a weight argument to see that $v_\nu^{i_0}$ must have the given form for some scalars $b(i_0, \ldots, i_r)$.  Then show that the definition of $b(i_0, \ldots, i_r)$ does give a \hwv\ by proving directly that $e_{k-1, k} v_\nu^{i_0}$ is zero for $1 < k \le m$.  This computation is also delicate.  For example, at $k > i_r$ three terms contribute multiples of
\begin{equation*}
   e_k \ot \bigl( e_{i_0, i_1} \cdots e_{i_{s-1}, k} \, e_{k-1, i_{s+1}} \cdots e_{i_{r-1}, i_r} v_\nu \bigr)
\end{equation*}
to $e_{k-1, k} v_\nu^{i_0}$.  Two of them correspond to the coefficients $b(i_0, \ldots, i_r)$ with $i_s = k-1$ and $i_s = k$, and the third corresponds to
\begin{equation*}
   b(i_0, \ldots, i_{s-1}, k, k-1, i_{s+1}, \cdots, i_r).
\end{equation*}
The reader may check that the multiples they contribute sum to zero.

Similarly, at $k = i_r$ multiples of $e_{k-1} \ot \bigl( e_{i_0, i_1} e_{i_1, i_2} \cdots e_{i_{r-1}, k} v_\nu \bigr)$ are contributed by three terms.  Again, the three multiples sum to zero.  $\Box$

\medbreak
Recall from Section~\ref{SMs} the \hwv\ $v_\delta^K$ in $\bC[\xi]_{(-k)} \ot L_\dlm(\delta)$ for $K \in \kappa(\delta, k)$.  Note that for any $d \le K_i$, the multi-index $K - d e_i$ is in $\kappa(\delta, k - d)$.

\begin{prop} \label{vdKd}
For $K \in \kappa(\delta, k)$ and $d \le K_i$, \th $L_\dlm(v_\delta^{K - d e_i}) \subseteq \Div^d \bigl( \F(v_\delta^K) \bigr)$.
\end{prop}

\meno {\it Proof.\/}
By Lemma~\ref{Divd}, $\Div^d \bigl( \bC[x]_{(d)} \ot L_\dlm(v_\delta^K) \bigr)$ is $\bC[\partial_\xi]_{(d)}$ applied to $L_\dlm(v_\delta^K)$, which is the same as $\bC[\partial_\xi]_{(1)}$ applied to $L_\dlm(v_\delta^K)$ $d$~times.  Thus it suffices to prove that $\bC[\partial_\xi]_{(1)}$ applied to $L_\dlm(v_\delta^K)$ contains $L_\dlm(v_\delta^{K - e_i})$ whenever $K_i > 0$.  In fact we need only prove that it contains $v_\delta^{K - e_i}$, because $\Div$ is $\dl_m$-invariant.

Consider $\bC[\partial_\xi]_{(1)} \ot L_\dlm(v_\delta^K)$.  By standard results on minuscule \r s, it contains a unique line of \hwv s of weight $\delta + \lambda_{K - e_{i_0}}$ for all $i_0$ with $K_{i_0} > 0$.  Applying Lemma~\ref{HWV} and recalling the identifications of $\dgl_m$ with $\dlm$ and $\bC^m$ with $\bC[\partial_\xi]_{(1)}$, we find that this line is spanned by
\begin{equation*}
   v_\delta^{K, i_0} := \sum_{\substack{{0 \le r < i_0} \\[1pt] {1 \le i_r < \cdots < i_1 < i_0}}} \,
   b(i_0, \ldots, i_r) \, \partial_{\xi_{i_r}} \ot \bigl( (x_{i_0} \partial_{x_{i_1}}) (x_{i_1} \partial_{x_{i_2}})
   \cdots (x_{i_{r-1}} \partial_{x_{i_r}}) v_\delta^K \bigr).
\end{equation*}
Here the factors $x_{i_{s-1}} \partial_{x_{i_s}}$ act on $v_\delta^K$ via the $\dlm$-action on $\bC[\xi]_{(-k)} \ot L_\dlm(\delta)$, and the scalars $b(i_0, \ldots, i_r)$ are defined using $\nu = \delta + \lambda_K$.  We have $\nu_{i_0 - 1} > \nu_{i_0}$ because $K_{i_0} > 0$.

We claim that the coefficient of $v_\delta$ in $\prod_{s=1}^r (x_{i_{s-1}} \partial_{x_{i_s}}) v_\delta^K$ is
\begin{equation*}
   (-1)^r \, (\xi_{i_1} \partial_{\xi_{i_0}}) (\xi_{i_2} \partial_{\xi_{i_1}}) \cdots (\xi_{i_r} \partial_{\xi_{i_{r-1}}}) \, \xi^K
   = (-1)^r \, K_{i_0} K_{i_1} \cdots K_{i_{r-1}} \, \xi^{K + e_{i_r} - e_{i_0}}.
\end{equation*}
To verify this, recall from Lemma~\ref{vdeltaK} that the coefficient of $v_\delta$ in $v_\delta^K$ is $\xi^K$.  Each operator $x_{i_{s-1}} \partial_{x_{i_s}}$ is of negative weight, so the only $\bC[\xi]$-multiple of $v_\delta$ arising as a summand of $\prod_{s=1}^r (x_{i_{s-1}} \partial_{x_{i_s}}) v_\delta^K$ occurs when every $x_{i_{s-1}} \partial_{x_{i_s}}$ acts on the factor $\xi^K$ of the lead term $\xi^K \ot v_\delta$.  Since $\phi_\xi(x_{i_{s-1}} \partial_{x_{i_s}}) = -\xi_{i_s} \partial_{\xi_{i_{s-1}}}$, the claim follows.

Now consider the map from $\bC[\partial_\xi]_{(1)} \ot \bC[\xi]_{(-k)} \ot L_\dlm(\delta)$ to $\bC[\xi]_{(1-k)} \ot L_\dlm(\delta)$ given by applying the first factor to the second.  This ``evaluation map'' is $\dlm$-covariant.  It must send $v_\delta^{K, i_0}$ to a multiple of $v_\delta^{K - e_{i_0}}$, because by Lemma~\ref{vdeltaK}, $v_\delta^{K - e_{i_0}}$ is the unique \hwv\ of its weight on the right side.  The proof of the proposition will be complete if we show that this multiple is non-zero.

By the last two displayed equations, the coefficient of $v_\delta$ in the image of $v_\delta^{K, i_0}$ under the evaluation map is
\begin{equation*}
   \sum_{\substack{{0 \le r < i_0} \\[1pt] {1 \le i_r < \cdots < i_1 < i_0}}} \,
   (-1)^r \, b(i_0, \ldots, i_r) \, K_{i_0} K_{i_1} \cdots K_{i_{r-1}} (K_{i_r} + 1) \, \xi^{K - e_{i_0}}.
\end{equation*}
Now observe that $c_j(i_0)^{-1} \in 2 + \bN$ for all $j < i_0$, so $(-1)^r b(i_0, \ldots, i_r)$ is always positive!  The contribution at $r = 0$ is non-zero, so the entire sum is a positive multiple of $\xi^{K - e_{i_0}}$.  $\Box$

\subsection{Projective operators} \label{amaps}

The point of this section is to use Lemma~\ref{bTFMops} and Proposition~\ref{vdKd} to show that the projective operators between \tfm s may be viewed as powers of the divergence operator.  These operators may be classified using the fact that the \tfm s are dual to parabolic Verma modules (see Proposition~\ref{dual Vermas}).  The dimension formulas for the spaces of homomorphisms between parabolic Verma modules derived in \cite{BES88} (see also \cite{ES86}, Section~6) give the following result.

\begin{thm} \label{TFM a maps}
$\dim \Hom_\dam \bigl( \F(\lambda), \F(\lambda') \bigr)$ is~0 unless $\chi_\lambda = \chi_{\lambda'}$ and $\lambda - \lambda' = d(\ep_0 - \ep_i)$ for some $d \in \bN$, when it is~1.  In this case, if\/ $d > 0$ then $\lambda = \mu[i-1]$ and $\lambda' = \mu[i]$ for some dominant integral $\mu$.
\end{thm}

Recall that $\S(\delta) = \bigoplus_{k=0}^\infty \S^k(\delta)$, and by~(\ref{decomp of Sk}), $\S^k(\delta) = \bigoplus_{\kappa(\delta, k)} \F(v_\delta^K)$.  Define projection operators $P_k$ and $P_K$ compatible with these $\vrm$-decompositions:
\begin{equation} \label{projections}
   P_k: \S(\delta) \twoheadrightarrow \S^k(\delta), \qquad
   P_K: \S(\delta) \twoheadrightarrow \F(v_\delta^K).
\end{equation}
These projections are related by the identity $P_k = \sum_{K \in \kappa(\delta, k)} P_K$.

\begin{prop} \label{Div powers}
For $K \in \kappa(\delta, k)$ and\/ $0 \le d \le K_i$,
\begin{equation*}
   \Hom_{\db_m} \bigl( \F(v_\delta^K), \F(v_\delta^{K - d e_i}) \bigr)
   = \bC P_{K - d e_i} \circ \Div^d \not= 0.
\end{equation*}
This map is\/ $\da_m$-covariant \iff\/ $\chi_{\delta + \lambda_K} = \chi_{\delta + \lambda_{K - d e_i}}$, \ie\ $\delta$, $K$, and $i$ satisfy~(\ref{RIC conditions}) and $K - d e_i$ is the partner of $K$.
\end{prop}

\meno {\it Proof.\/}
The second sentence follows from the first, Proposition~\ref{RIC Prop}, and Theorem~\ref{TFM a maps}.  Observe that the first sentence generalizes the example in Remark~\ref{Sec 7.1 Remark 1}.  Since $\Div$ is a $\dbm$-map and $P_{K - d e_i}\big|_{\S^{k-d}(\delta)}$ is a $\vrm$-map, $P_{K - d e_i} \circ \Div^d$ is a $\dbm$-map, and it is non-zero by Proposition~\ref{vdKd}.  It will suffice to prove that it is up to a scalar unique.

By Lemmas~\ref{vdeltaK} and~\ref{bTFMops}, $\Hom_{\db_m} \bigl( \F(v_\delta^K), \F(v_\delta^{K - d e_i}) \bigr)$ is isomorphic to
\begin{equation*}
   \Bigl[ \bC[\px]_{(-d)} \ot L_\dlm(\delta + \lambda_K)^*
   \ot L_\dlm(\delta + \lambda_{K - d e_i}) \Bigr]^{\dl_m}.
\end{equation*}
By the Parthsarathy-Ranga Rao-Varadarajan lemma, $L_\dlm(d \ep_1 - d \ep_0)$ occurs with multiplicity~1 in $L_\dlm(\delta + \lambda_K)^* \ot L_\dlm(\delta + \lambda_{K - d e_i})$, as its smallest submodule.  This module is dual to $\bC[\px]_{(-d)} \simeq L_\dlm(d \ep_0 - d \ep_m)$, so the proof is complete.  $\Box$

\section{The Jordan form of the Casimir operator}  \label{JORDAN} 

\def\gd{{\gamma, \delta}}

In this section we determine those \ic s $\chi_\mu$ of $\dam$ such that the generalized $\chi_\mu$-submodule $\D_\gamma(\delta)^{(\mu)}$ is not equal to the $\chi_\mu$-submodule $\D_\gamma(\delta)^\mu$.  By Proposition~\ref{folklore}, these are precisely the \ic s for which the Casimir operator is not semisimple.

By Proposition~\ref{RIC Prop}, the symbol module $\S(\delta)$ of $\D_\gamma(\delta)$ has at most two submodules $\F(v_\delta^K)$ with any given \ic.  Therefore
\begin{equation*}
   \D_\gamma(\delta)^{(\mu)} = \D_\gamma(\delta)^{(\mu: 2)}
\end{equation*}
for all $\mu$.  It follows from Propositions~\ref{par-injectives} and~\ref{Shelton} that if $\D_\gamma(\delta)^{(\mu)} \not= \D_\gamma(\delta)^\mu$, then we may assume without loss of generality that $\mu$ is dominant integral and $\D_\gamma(\delta)^{(\mu)} = I(\mu[i])$ for some $1 \le i \le m$.

For convenience we abbreviate the $\vrm$-actions on $\D_\gamma(\delta)$ and $\S(\delta)$: set
\begin{equation*}
   \Lie_\gd := \Lie_{\bC_\gamma, \bC_\gamma \ot L_\dlm(\delta)}, \qquad
   \Lie^\S_\delta := \Lie^\S_{\bC_\gamma, \bC_\gamma \ot L_\dlm(\delta)}.
\end{equation*}
In order to take advantage of Proposition~\ref{LpppOmegaa}, instead of working with $\Lie_\gd$ we will work with the isomorphic action $\Lie_\gd^{\NS}$ on $\S(\delta)$.

Recall from~(\ref{projections}) the projections $P_k$ and $P_K$, and set $\S_k(\delta) := \bigoplus_{j=0}^k \S^j(\delta)$.  By Lemma~\ref{LieVV}, $\S_k(\delta)$ is invariant under $\Lie_\gd^{\NS}$, and the restrictions
\begin{equation*}
   P_k\big|_{\S_k(\delta)}: \S_k(\delta) \twoheadrightarrow \S^k(\delta) \qquad
   P_K\big|_{\S_k(\delta)}: \S_k(\delta) \twoheadrightarrow \F(v_\delta^K)
\end{equation*}
are $\vrm$ intertwining maps from $\Lie^{\NS}_\gd$ to $\Lie^\S_\delta$.

{\em In this section, when we write $\S(\delta)$ or $\S_k(\delta)$ without specifying an action, it is understood that the action is\/ $\Lie^{\NS}_\gd$.}

\begin{lemma} \label{no partner}
Suppose that $K \in \kappa(\delta)$ does not satisfy~(\ref{RIC conditions}), and set $k := |K|$.
\begin{enumerate}

\item[(i)]
$\S(\delta)^{(\delta + \lambda_K)} = \S_k(\delta)^{\delta + \lambda_K}$, and $\S_{k - 1}(\delta)^{(\delta + \lambda_K)} = 0$.

\smallbreak \item[(ii)]
$P_k$ and $P_K$ map\/ $\S_k(\delta)^{\delta + \lambda_K}$ bijectively to\/ $\F(v_\delta^K)$.

\end{enumerate}
\end{lemma}

\meno {\it Proof.\/}
For all weights $\nu$, the projection $P_j$ carries $\S_j(\delta)^{(\nu)}$ onto $\S^j(\delta)^{(\nu)}$ with kernel $\S_{j-1}(\delta)^{(\nu)}$.  The space $\S^j(\delta)^{(\mu)}$ is~$0$ unless $j = k$, when it is $\F(v_\delta^K)$.  The lemma follows by induction on order.  $\Box$

\medbreak
For the remainder of this section, fix $\delta$, $K \in \kappa(\delta)$, and~$i$ satisfying~(\ref{RIC conditions}).  Let $K'$ be the partner of $K$, assume $|K| > |K'|$, and set $k := |K|$ and $k' := |K'|$.  Let $\mu$ be the dominant integral weight such that $\delta + \lambda_K = \mu[i-1]$ and $\delta + \lambda_{K'} = \mu[i]$, as in Proposition~\ref{RIC Prop}(iii).  

By Proposition~\ref{RIC Prop}, $K$ and $K'$ are the only multi-indices $J \in \kappa(\delta)$ such that $\chi_{\delta + \lambda_J} = \chi_\mu$.  Therefore Proposition~\ref{Decomp of Sk} gives
\begin{equation*}
   \bigl( \S(\delta),\, \Lie^\S_\delta \bigr)^{(\mu)} =
   \bigl( \S(\delta),\, \Lie^\S_\delta \bigr)^{\mu} =
   \F(v_\delta^K) \oplus \F(v_\delta^{K'}).
\end{equation*}

\begin{lemma} \label{D IC espaces}
\begin{enumerate}
\item[(i)]
$\S(\delta)^{(\mu)} = \S_k(\delta)^{(\mu: 2)}$.

\smallbreak \item[(ii)]
$P_k$ and $P_K$ map\/ $\S_k(\delta)^{(\mu: 2)}$ onto\/ $\F(v_\delta^K)$ with kernel\/ $\S_{k'}(\delta)^\mu$.

\smallbreak \item[(iii)]
$\S_{k-1}(\delta)^{(\mu)} = \S_{k'}(\delta)^\mu$, and $\S_{k' - 1}(\delta)^{(\mu)} = 0$.

\smallbreak \item[(iv)]
$P_{k'}$ and $P_{K'}$ map\/ $\S_{k'}(\delta)^\mu$ bijectively to\/ $\F(v_\delta^{K'})$.

\end{enumerate}
\end{lemma}

\meno {\it Proof.\/}
Induct on order as in the proof of Lemma~\ref{no partner}.  Use the fact that $\S^j(\delta)^{(\mu)}$ is~$0$ unless $j$ is $k$ or $k'$, when it is $\F(v_\delta^K)$ or $\F(v_\delta^{K'})$, respectively.  $\Box$

\medbreak
The next theorem is a rephrasing of Theorem~\ref{AM Decomp}(ii).  Its proof occupies the remainder of this section.

\begin{thm} \label{baby main theorem}
\begin{enumerate}
\item[(i)]
For $k + (m + 1) \gamma \not\in \{1, \ldots, k - k' \}$,\/ $\S(\delta)^{(\mu)} \simeq I(\mu[i])$.

\smallbreak \item[(ii)]
For $k + (m + 1) \gamma \in \{1, \ldots, k - k' \}$,\/ $\S(\delta)^{(\mu)} \simeq \F(\mu[i-1]) \oplus \F(\mu[i])$.

\end{enumerate}
\end{thm}

The first step in the proof is to recall the {\em tree-like subspaces\/} $\T(v_\delta^J) \subset \S(\delta)$ defined in \cite{MR07}.  Put the standard partial order on $\bN^m$:
\begin{equation*}
   J' \le J \mbox{\rm\ \ if\ \ } J'_i \le J_i \mbox{\rm\ \ for\ \ } 1 \le i \le m.  
\end{equation*}
Note that if $J \in \kappa(\delta)$, then $J' \in \kappa(\delta)$ for all $J' \le J$.  Set
\begin{equation*}
   \T(v_\delta^J) := \bigoplus_{J' \le J} \F(v_\delta^{J'}).
\end{equation*}

\begin{lemma} \label{MaRa} \cite{MR07}
$\T(v_\delta^J)$ is invariant under\/ $\Lie^\S_\delta(\vrm)$ and\/ $\Lie^{\NS}_\gd(\dam)$.
\end{lemma}

\meno {\it Proof.\/}
The $\Lie^\S_\delta(\vrm)$-invariance follows from that of the summands $\F(v_\delta^{J'})$.  For the $\Lie^{\NS}_\gd(\dam)$-invariance, use~(\ref{LieVVa}) to obtain
\begin{equation} \label{JNT}
   \Lie^{\NS}_\gd (x_i \E_x) - \Lie^\S_\delta (x_i \E_x) = -\bigl( \E_\xi + (m + 1) \gamma \bigr) \partial_{\xi_i}.
\end{equation}
Recall that $\F(v_\delta^{J'}) = \bC[x] \ot L_\dlm(v_\delta^{J'})$, and reduce to proving
\begin{equation} \label{image of Div}
   \partial_{\xi_i} L_\dlm(v_\delta^{J'}) \subseteq \bigoplus_{j:\th J'_j > 0} L_\dlm(v_\delta^{J' - e_j}).
\end{equation}

As we saw in the proof of Proposition~\ref{vdKd}, the \hw s of the $\dlm$-\irr\ summands of $\bC[\partial_\xi]_{(1)} \ot L_\dlm(v_\delta^{J'})$ are all of the form $\delta + \lambda_{J'} - \ep_0 + \ep_j$.  Evaluating the $\xi$-derivatives and matching \hw s with~(\ref{decomp of Sk}) completes the proof.  $\Box$

\medbreak
For any $\nu \in \dh_m^*$, define $P^{(\nu)}$ to be the projection of $\S(\delta)$ to its generalized $\chi_\nu$-submodule $\S(\delta)^{(\nu)}$ along the sum $\bigoplus_{\nu' \not= \nu} \S(\delta)^{(\nu')}$ of its other generalized $\chi_{\nu'}$-submodules.  The action of $P^{(\delta + \lambda_J)}$ on $\F(v_\delta^J)$ is~$1$ at the symbol level:
\begin{equation} \label{IC projections}
   P_{|J|} \circ P^{(\delta + \lambda_J)}\big|_{\F(v_\delta^J)} =
   P_{J} \circ P^{(\delta + \lambda_J)}\big|_{\F(v_\delta^J)} = 1.
\end{equation}

It is elementary that $P^{(\nu)} = \Lie^{\NS}_\gd(\Upsilon)$ for some $\Upsilon \in \dZ(\dam)$.  Coupled with Lemma~\ref{MaRa} this yields the next lemma, a sharpening of Lemma~\ref{D IC espaces}.

\begin{lemma} \label{Tree IC espaces}
\begin{enumerate}
\item[(i)]
$\S(\delta)^{(\mu)} = \T(v_\delta^K)^{(\mu: 2)} = P^{(\mu)} \bigl( \F(v_\delta^K) \oplus \F(v_\delta^{K'}) \bigr)$.

\smallbreak \item[(ii)]
$\S_{k'}(\delta)^\mu = \T(v_\delta^{K'})^\mu = P^{(\mu)} \bigl(\F(v_\delta^{K'}) \bigr)$.

\end{enumerate}
\end{lemma}

\begin{remark} \label{Sec 8 Remark}
Under the action $\Lie_\gd^\S$, $\F(v_\delta^{K'})$ is the only submodule of\/ $\T(v_\delta^K)$ with the same \ic\ as $\F(v_\delta^K)$.  However, there can exist other submodules $\F(v_\delta^J)$ with the same Casimir eigenvalue.  Using (\ref{chilambdaOmega}), one sees that an example is provided by
\begin{equation*}
   m=2, \quad \delta = (-4, 4, 0), \quad K = (3, 2), \quad K' = (2, 2), \quad J = (1, 1).
\end{equation*}
\end{remark}

We now need three technical results describing the interaction between projections and central operators, in particular, the Casimir operator.  It will be convenient to use the convention that $P_J = 0$ if $J \not\in \kappa(\delta, |J|)$.  Define
\begin{equation*}
   \nu_c := \delta + \lambda_{K - c e_i}, \qquad
   \omega_c := \chi_{\nu_c} (\Omega_\dam), \qquad d := k - k'.
\end{equation*}

\begin{lemma} \label{Projections & Div}
Assume $0 \le c \le d$ and $\Theta \in \dZ(\dam)$.  Abbreviate
\begin{equation*}
   \Theta^{\NS} := \Lie^{\NS}_\gd(\Theta),\ \
   \Omega^{\NS}_c := \Lie^{\NS}_\gd(\Omega_\dam - \omega_c),\ \
   q_c := 2 \bigl( k - c -1 + (m + 1) \gamma \bigr).
\end{equation*}

\begin{enumerate}
\item[(i)]
$P_{K'} \circ \Theta^{\NS} \circ P_J = 0$ if $K' \not\le J$, in particular, if $J \le K$ and $K - J \not\in \bN e_i$.

\medbreak \item[(ii)]
$P_{K'} \circ \Theta^{\NS} \circ \bigl( \prod_{a=c}^d \Omega^{\NS}_a \bigr) \circ P_{K - c e_i} = 0$.

\medbreak \item[(iii)]
$P_{K'} \circ \bigl( \prod_{a=0}^{d-1} \Omega^{\NS}_a \bigr) \circ P_K = \bigl( \prod_{a=0}^{d-1} q_a \bigr) P_{K'} \circ \Div^{d} \circ P_K$.

\end{enumerate}
\end{lemma}

\meno {\it Proof.\/}
By Lemma~\ref{MaRa}, the image of $\Theta^{\NS} \circ P_J$ is in $\T(v_\delta^J)$.  For $K' \not\le J$, $\F(v_\delta^{K'})$ is not a summand of $\T(v_\delta^J)$, proving~(i).

In order to proceed we must prove a series of equations: for $c < b \le d+1$,
\begin{align}
   & \Omega^{\NS}_c \circ P_{K - c e_i}
   = q_c \Div \circ P_{K - c e_i}
   = q_c \bigl( \sum_{j = 1}^m P_{K - c e_i - e_j} \bigr) \circ \Div \circ P_{K - c e_i},
   \label{ii} \\
   & P_{K'} \circ \Theta^{\NS} \circ \Omega^{\NS}_c \circ P_{K - c e_i}
   = q_c P_{K'} \circ \Theta^{\NS} \circ P_{K - (c+1) e_i} \circ \Div \circ P_{K - c e_i},
   \label{iii} \\
   & P_{K'} \circ \Theta^{\NS} \circ
   \bigl( \prod_{a=c}^{b-1} \Omega^{\NS}_a \bigr) \circ P_{K - c e_i}
   = \bigl( \prod_{a=c}^{b-1} q_a \bigr) P_{K'}
   \circ \Theta^{\NS} \circ \Div^{b-c} \circ P_{K - c e_i}.
   \label{iv}
\end{align}

The first equality in~(\ref{ii}) is immediate from Lemma~\ref{LpppOmegaa}(iii) and the definition of $\omega_c$.  For the second, note that the argument used to prove~(\ref{image of Div}) shows that $\Div$ maps $\F(v_\delta^{K - ce_i})$ to $\bigoplus_j \F(v_\delta^{K - ce_i - e_j})$.

For~(\ref{iii}), apply~(\ref{ii}) and then use~(i) to see that the terms other than $P_{K - (c+1) e_i}$ in the sum $\sum_{j = 1}^m P_{K - c e_i - e_j}$ contribute nothing.

For~(\ref{iv}), apply~(\ref{iii}) repeatedly to obtain $P_{K'} \circ \Theta^{\NS} \circ \bigl( \prod_{a=c}^{b-1} \Omega^{\NS}_a \bigr) \circ P_{K - c e_i} = $
\begin{equation*}
   \Bigl( \prod_{a=c}^{b-1} q_a \Bigr) P_{K'} \circ \Theta^{\NS} \circ P_{K - b e_i} \circ
   \Div \circ P_{K - (b-1) e_i} \circ \Div \circ \cdots \circ P_{K - (c+1) e_i} \circ \Div \circ P_{K - c e_i}.
\end{equation*}
Then use~(i) and the second equality of~(\ref{ii}) to prove that this operator remains the same if the internal projections are dropped.

Now we can prove~(ii): apply~(\ref{iv}) with $b = d+1$, and note that $\Theta^{\NS} \circ \Div^{d-c+1}$ maps $\F(v_\delta^{K - ce_i})$ to $\S_{k'-1}(\delta)$, which is annihilated by $P_{K'}$.

Finally, (iii) is~(\ref{iv}) with $c=0$, $b=d$, and $\Theta^{\NS} = 1$.  $\Box$

\begin{lemma} \label{bar}
Let\/ $\Theta$ be any element of\/ $\dZ(\dam)$ such that\/ $\chi_\mu(\Theta) = 0$.  Then there exists a unique\/ $\dam$-covariant map\/ $\o\Theta: \F(v_\delta^K) \to \F(v_\delta^{K'})$ such that
\begin{equation*}
   \o\Theta \circ P_K \big|_{\S_k(\delta)^{(\mu)}} =
   P_{K'} \circ \Lie^{\NS}_\gd (\Theta) \big|_{\S_k(\delta)^{(\mu)}}.
\end{equation*}
\end{lemma}

\meno {\it Proof.\/}
This is a commutative diagram lemma.  Since $\chi_\mu(\Theta) = 0$, $\Lie^\S_\delta(\Theta)$ annihilates $\F(v_\delta^K)$ and $\F(v_\delta^{K'})$.  Hence by Lemma~\ref{D IC espaces}, $\Lie^{\NS}_\gd (\Theta)$ maps $\S_k(\delta)^{(\mu)}$ to $\S_{k'}(\delta)^\mu$ and annihilates $\S_{k'}(\delta)^\mu$, and its restriction to $\S_k(\delta)^{(\mu)}$ factors through $P_K$.  The rest follows easily.  $\Box$

\medbreak
In order to apply Lemma~\ref{bar}, we must choose $\Theta$ so that we can compute $\o\Theta$ explicitly.  The lemma applies to any $\dZ(\dam)$-multiple of $\Omega_\dam - \omega_0$, but if we take for example $\Theta$ to be $\Omega_\dam - \omega_0$, we have no way to compute $\o\Theta$.  Keeping in mind that $\omega_0 = \omega_d$, define
\begin{equation*}
   \Delta := \bigl( \oh \bigr)^d \prod_{c=0}^{d-1} (\Omega_\dam - \omega_c)
   = \bigl( \oh \bigr)^d \prod_{c=1}^d (\Omega_\dam - \omega_c).
\end{equation*}

\begin{prop} \label{Delta}
$\o\Delta = \Bigl( \prod_{c=1}^d \bigl[ k - c + (m + 1) \gamma \bigr] \Bigr) P_{K'} \circ \Div^d$.
\end{prop}

\meno {\it Proof.\/}
Fix $v \in \F(v_\delta^K)$ and write $\t v$ for $P^{(\nu_0)} v$.  By~(\ref{IC projections}), $P_K \t v = v$, so by Lemma~\ref{bar}, $\o\Delta (v) = P_{K'} \circ \Lie^{\NS}_\gd(\Delta) (\t v)$.  By Lemma~\ref{Tree IC espaces}, $\t v = \sum_{J \le K} P_J \t v$, and by Lemma~\ref{Projections & Div}(i) and~(ii) (note that $\Omega^{\NS}_0 = \Omega^{\NS}_d$), $P_{K'} \circ \Lie^{\NS}_\gd(\Delta) (P_J \t v) = 0$ for $J < K$.  Therefore $\o\Delta (v) = P_{K'} \circ \Lie^{\NS}_\gd(\Delta) (v)$.  Lemma~\ref{Projections & Div}(iii) now completes the proof.  $\Box$

\begin{remark}
This proposition gives an alternate proof of the fact that the map $P_{K'} \circ \Div^d: \F(v^K_\delta) \to \F(v^{K'}_\delta)$ is $\dam$-covariant, independent of Proposition~\ref{Div powers}.
\end{remark}

\meno {\it Proof of Theorem~\ref{baby main theorem}.\/}
By Proposition~\ref{vdKd}, $P_{K'} \circ \Div^d \bigl( \F(v_\delta^K) \bigr) \not= 0$.  Therefore by Proposition~\ref{Delta}, if $k + (m + 1) \gamma \not\in \{1, \ldots, d \}$, then $\o\Delta \not= 0$.  In this case $\S(\delta)^{(\nu_0)}$ is not annihilated by $\Omega_\dam - \omega_0$, so $\Omega_\dam$ does not act semisimply on it.  Hence by Proposition~\ref{Shelton} and Lemma~\ref{D IC espaces}, $\S(\delta)^{(\nu_0)} \simeq I(\mu[i])$.  This proves~(i).

In~(ii), $\o\Delta = 0$, so an elementary argument shows that $\Omega_\dam$ acts by $\omega_0$ on $\S(\delta)^{(\nu_0)}$.  Hence $\S(\delta)^{(\nu_0)} \not= I(\mu[i])$ by~(iv), so by Proposition~\ref{Shelton} it splits as desired.  (For an alternate proof of~(ii), one can use the fact that by~(\ref{JNT}), if\/ $-(m+1) \gamma \in \bN$, then $\bigoplus_{j = 1 - (m+1) \gamma}^\infty \S^j(\delta)$ is invariant under $\Lie_\gd^{\NS}|_\dam$.)  $\Box$

\section{Proofs}  \label{PROOFS} 

Here we tie together the results of Sections~\ref{SC}-\ref{JORDAN} to prove the results of Section~\ref{MRs}.  We use the notation of Remark~\ref{tdelta}.  Recall that Propositions~\ref{Decomp of Sk} and~\ref{RIC Prop} were proven in Section~\ref{SMs}.

\meno {\it Proof of Theorem~\ref{AM Decomp}.\/}
Keeping in mind that the actions $\Lie_\gd$ on $\D_\gamma(\delta)$ and $\Lie_\gd^{\NS}$ on $\S(\delta)$ are isomorphic by construction, (i) is simply a restatement of Lemma~\ref{no partner}, and~(ii)(a)-(c) restate Lemma~\ref{D IC espaces} and Theorem~\ref{baby main theorem}.  For~(ii)(d), use $\lambda_K - \lambda_{K'} = \mu[i-1] - \mu[i]$ coupled with~(\ref{mimommi}).  $\Box$

\meno {\it Proof of Corollary~\ref{main cor}.\/}
If $\delta$ is either dominant or singular, $K = 0$ fails~(\ref{RIC conditions}), so it has no partner.  Hence Theorem~\ref{AM Decomp}(i) gives~(i).

For~(ii), use~(\ref{mimommi}) to verify that $K = 0$ satisfies~(\ref{RIC conditions}) and has partner $K' = (\mu_{i-1} - \mu_i + 1) e_i$.  Theorem~\ref{AM Decomp}(ii) then completes the proof.  $\Box$

\meno {\it Proof of Proposition~\ref{Splitting Prop}.\/}
If~(ii) fails, (\ref{GIC Decomp}) gives the unique $\dam$-splitting of~(\ref{exact seq}):
\begin{equation*}
   \D^k_\gamma(\delta) = \D^{k-1}_\gamma(\delta) \oplus
   \bigoplus_{K \in \kappa(\delta, k)} \D^k_\gamma(\delta)^{\delta + \lambda_K}.
\end{equation*}
If~(ii) holds, then either~(a) or~(b) of Theorem~\ref{AM Decomp}(ii) holds.  If~(a), then~(\ref{exact seq}) does not split because~(\ref{I no split}) does not split.  If~(b), then~(\ref{exact seq}) splits non-uniquely because by Theorem~\ref{TFM a maps} there are non-trivial $\dam$-maps from $\F(\mu[i-1])$ to $\F(\mu[i])$.  Thus~(i) and~(ii) are equivalent.  Proposition~\ref{RIC Prop} gives the equivalence of~(ii) and~(iii).  

If~(iii) holds, then Proposition~\ref{RIC Prop} entails $\t\delta_i \ge k > |K'| > \t\delta_{i+1}$ and $2k > k + K'_i = \t\delta_i$, where $i = i(\delta, k)$.  These imply~(iv).  Conversely, if~(iv) holds, set $K := K_i e_i + (k - K_i) e_m$, where $K_i := \min \{k, \delta_i - \delta_{i+1} \}$.  Check that $K$ is in $\kappa(\delta, k)$ and satisfies~(\ref{RIC conditions}), and its partner $K' = (\t\delta_i - k) e_i + (k - K_i) e_m$ satisfies~(iii).    $\Box$

\begin{remark}
Our approach to Proposition~\ref{Splitting Prop} is independent of\/ \cite{Mi12} and gives an alternate proof of his Theorem~\ref{JPM Thm} for $\D_\gamma(\delta)$.  Indeed, by Proposition~\ref{RIC Prop} and Theorem~\ref{TFM a maps} there exist non-trivial $\dam$-maps between distinct symbol modules $\S^k(\delta)$ and $\S^{k'}(\delta)$ \iff\/ $\S(\delta)$ has repeated \ic s.  By Proposition~\ref{Resonance Prop}, this occurs \iff\/ $\D_\gamma(\delta)$ is resonant.
\end{remark}

\meno {\it Proof of Proposition~\ref{Resonance Prop}.\/}
$\D_\gamma(\delta)$ is resonant \iff~(\ref{exact seq}) fails to have a unique $\dam$-splitting for some $k$.  By Proposition~\ref{Splitting Prop}, this occurs \iff\ $\S(\delta)$ has repeated \ic s.  Thus~(i) and~(ii) are equivalent.

Suppose that $\S(\delta)$ has repeated \ic s $\chi_{\delta + \lambda_K}$ and $\chi_{\delta + \lambda_{K'}}$, with $|K| > |K'|$.  Writing $i$ for $i(\delta, |K|)$, we have $\t\delta_i \ge |K| > |K'| > \t\delta_{i+1}$, which gives both $\delta_i > \delta_{i+1}$ and $\t\delta_i > 0$.  It follows that~(ii) implies~(iii).  Conversely, if~(iii) holds, let $i$ be maximal such that $\delta_i = \delta_1$.  Then $k = \t\delta_i$ satisfies~(\ref{Splitting conditions}), so~(ii) holds.  (In fact, here $e_i + (\t\delta_i - 1) e_m$ has partner $(\t\delta_i - 1) e_m$ if $i < m$, and $0$ if $i = m$.)  $\Box$

\meno {\it Proof of Theorem~\ref{Multi PQ Thm}.\/}
Because~(\ref{exact seq}) does not have a unique splitting, Propositions~\ref{RIC Prop} and~\ref{Splitting Prop} imply that~(\ref{RIC conditions}) holds.  Writing~$i$ for $i(\delta, k)$, they also imply that $2k > \t\delta_i$ and $\t\delta_i \ge k > \t\delta_{i+1} + 1$, and that there is at least one multi-index $K \in \kappa(\delta, k)$ having a partner $K'$ with $|K'| < k$.

For all such $K$, $\D^k_\gamma(\delta)^{(\delta + \lambda_K)}$ is equivalent to either $I(\delta + \lambda_{K'})$ or $\F(\delta + \lambda_K) \oplus \F(\delta + \lambda_{K'})$.  If for any $K$ it is equivalent to $I(\delta + \lambda_{K'})$, then~(\ref{exact seq}) has no $\dam$-splitting because~(\ref{I no split}) has none.  Otherwise~(\ref{exact seq}) has a continuous family of splittings, because there are non-trivial $\dam$-maps from $\F(\delta + \lambda_K)$ to $\F(\delta + \lambda_{K'})$.  This proves~(i).

For~(ii) and~(iii) we apply Theorem~\ref{AM Decomp}(ii).  Suppose $i < m$.  If $k + (m+1) \gamma = 1$, then for all $K$ as above, $\D^k_\gamma(\delta)^{(\delta + \lambda_K)}$ splits.  Conversely, check that
\begin{equation*}
   K := (\t\delta_i - k + 1) e_i + (2k - \t\delta_i - 1) e_m, \quad
   K' := (\t\delta_i - k) e_i + (2k - \t\delta_i - 1) e_m
\end{equation*}
are partners with $|K| = k$ and $|K'| = k-1$.  If $k + (m+1) \gamma \not= 1$, $\D^k_\gamma(\delta)^{(\delta + \lambda_K)}$ does not split for this $K$.  This proves~(ii).

Finally, suppose $i = m$, and write $d$ for $\max \{ 1, 2k + m + \delta_0 - \delta_1 \}$.  For $K$ and $K'$ as above, we know that $\D^k_\gamma(\delta)^{(\delta + \lambda_K)}$ splits \iff\ $k + (m+1) \gamma$ is in $\{ 1, \ldots, k - |K'| \}$.  Here $|K'| = \t\delta_m - K_m$ and $K_m = k - \sum_1^{m-1} K_i$.  Because $K_i \le \delta_i - \delta_{i+1}$, we have $K_m \ge k - \delta_1 + \delta_m$.  Thus
\begin{equation*}
   k - |K'| = k + K_m - \t\delta_m \ge 2k + m + \delta_0 - \delta_1.
\end{equation*}
Therefore, if $k + (m+1) \gamma$ is in $\{1, \ldots, d \}$, $\D^k_\gamma(\delta)^{(\delta + \lambda_K)}$ splits for all allowed $K$.

Conversely, suppose that $k + (m + 1) \gamma$ is not in $\{ 1, \ldots, d \}$.  For $d > 1$, take
\begin{equation*}
   K = (\delta_1 - \delta_2) e_1 + \cdots (\delta_{m-1} - \delta_m) e_{m-1} + (k - \delta_1 + \delta_m) e_m,
   \quad K' = K - d e_m.
\end{equation*}
Then $K$ and $K'$ are partners for which $\D^k_\gamma(\delta)^{(\delta + \lambda_K)}$ does not split.

For $d = 1$, construct partners $K$ and $K' = K - e_m$: take $K_m = \t\delta_m - k +1$ and choose $K_1, \ldots, K_{m-1}$ in any way satisfying $K_i \le \delta_i - \delta_{i+1}$ and $|K| = k$.  Then again, $\D^k_\gamma(\delta)^{(\delta + \lambda_K)}$ does not split.  This proves~(iii).  $\Box$

\meno {\bf Acknowledgments.}
We are grateful to V.~S.~Varadarajan for encouraging us to seek the injectives of $\O^{\dgl_m}(\dsl_{m+1})$ in the \dog\ modules of $\vrm$.  In addition, we thank V.~Mazorchuk for helpful discussions, E.~Mukhin for advice concerning the exposition of Section~\ref{bmaps}, and the referee for improving the accessibility of the paper.

\bibliographystyle{amsalpha}

\end{document}